wsp


\input amstex
\expandafter\ifx\csname mathdefs.tex\endcsname\relax
  \expandafter\gdef\csname mathdefs.tex\endcsname{}
\else \message{Hey!  Apparently you were trying to
  \string twice.   This does not make sense.} 
\errmessage{Please edit your file (probably \jobname.tex) and remove
any duplicate ``\string\input'' lines} \fi




\catcode`\X=12\catcode`\@=11

\def\n@wcount{\alloc@0\count\countdef\insc@unt}
\def\n@wwrite{\alloc@7\write\chardef\sixt@@n}
\def\n@wread{\alloc@6\read\chardef\sixt@@n}
\def\r@s@t{\relax}\def\v@idline{\par}\def\@mputate#1/{#1}
\def\l@c@l#1X{\firstpart.#1}\def\gl@b@l#1X{#1}\def\t@d@l#1X{{}}

\def\crossrefs#1{\ifx\all#1\let\tr@ce=\all\else\def\tr@ce{#1,}\fi
   \n@wwrite\cit@tionsout\openout\cit@tionsout=\jobname.cit 
   \write\cit@tionsout{\tr@ce}\expandafter\setfl@gs\tr@ce,}
\def\setfl@gs#1,{\def\@{#1}\ifx\@\empty\let\next=\relax
   \else\let\next=\setfl@gs\expandafter\xdef
   \csname#1tr@cetrue\endcsname{}\fi\next}
\def\m@ketag#1#2{\expandafter\n@wcount\csname#2tagno\endcsname
     \csname#2tagno\endcsname=0\let\tail=\all\xdef\all{\tail#2,}
   \ifx#1\l@c@l\let\tail=\r@s@t\xdef\r@s@t{\csname#2tagno\endcsname=0\tail}\fi
   \expandafter\gdef\csname#2cite\endcsname##1{\expandafter
     \ifx\csname#2tag##1\endcsname\relax?\else\csname#2tag##1\endcsname\fi
     \expandafter\ifx\csname#2tr@cetrue\endcsname\relax\else
     \write\cit@tionsout{#2tag ##1 cited on page \folio.}\fi}
   \expandafter\gdef\csname#2page\endcsname##1{\expandafter
     \ifx\csname#2page##1\endcsname\relax?\else\csname#2page##1\endcsname\fi
     \expandafter\ifx\csname#2tr@cetrue\endcsname\relax\else
     \write\cit@tionsout{#2tag ##1 cited on page \folio.}\fi}
   \expandafter\gdef\csname#2tag\endcsname##1{\expandafter
      \ifx\csname#2check##1\endcsname\relax
      \expandafter\xdef\csname#2check##1\endcsname{}%
      \else\immediate\write16{Warning: #2tag ##1 used more than once.}\fi
      \multit@g{#1}{#2}##1/X%
      \write\t@gsout{#2tag ##1 assigned number \csname#2tag##1\endcsname\space
      on page \number\count0.}%
   \csname#2tag##1\endcsname}}
\def\multit@g#1#2#3/#4X{\def\t@mp{#4}\ifx\t@mp\empty%
      \global\advance\csname#2tagno\endcsname by 1 
      \expandafter\xdef\csname#2tag#3\endcsname
      {#1\number\csname#2tagno\endcsnameX}%
   \else\expandafter\ifx\csname#2last#3\endcsname\relax
      \expandafter\n@wcount\csname#2last#3\endcsname
      \global\advance\csname#2tagno\endcsname by 1 
      \expandafter\xdef\csname#2tag#3\endcsname
      {#1\number\csname#2tagno\endcsnameX}
      \write\t@gsout{#2tag #3 assigned number \csname#2tag#3\endcsname\space
      on page \number\count0.}\fi
   \global\advance\csname#2last#3\endcsname by 1
   \def\t@mp{\expandafter\xdef\csname#2tag#3/}%
   \expandafter\t@mp\@mputate#4\endcsname
   {\csname#2tag#3\endcsname\lastpart{\csname#2last#3\endcsname}}\fi}
\def\t@gs#1{\def\all{}\m@ketag#1e\m@ketag#1s\m@ketag\t@d@l p
   \m@ketag\gl@b@l r \n@wread\t@gsin
   \openin\t@gsin=\jobname.tgs \re@der \closein\t@gsin
   \n@wwrite\t@gsout\openout\t@gsout=\jobname.tgs }
\outer\def\localtags{\t@gs\l@c@l}
\outer\def\globaltags{\t@gs\gl@b@l}
\outer\def\newlocaltag#1{\m@ketag\l@c@l{#1}}
\outer\def\newglobaltag#1{\m@ketag\gl@b@l{#1}}

\newif\ifpr@ 
\def\m@kecs #1tag #2 assigned number #3 on page #4.%
   {\expandafter\gdef\csname#1tag#2\endcsname{#3}
   \expandafter\gdef\csname#1page#2\endcsname{#4}
   \ifpr@\expandafter\xdef\csname#1check#2\endcsname{}\fi}
\def\re@der{\ifeof\t@gsin\let\next=\relax\else
   \read\t@gsin to\t@gline\ifx\t@gline\v@idline\else
   \expandafter\m@kecs \t@gline\fi\let \next=\re@der\fi\next}
\def\pretags#1{\pr@true\pret@gs#1,,}
\def\pret@gs#1,{\def\@{#1}\ifx\@\empty\let\n@xtfile=\relax
   \else\let\n@xtfile=\pret@gs \openin\t@gsin=#1.tgs \message{#1} \re@der 
   \closein\t@gsin\fi \n@xtfile}

\newcount\sectno\sectno=0\newcount\subsectno\subsectno=0
\newif\ifultr@local \def\ultralocal{\ultr@localtrue}
\def\firstpart{\number\sectno}
\def\lastpart#1{\ifcase#1 \or a\or b\or c\or d\or e\or f\or g\or h\or 
   i\or k\or l\or m\or n\or o\or p\or q\or r\or s\or t\or u\or v\or w\or 
   x\or y\or z \fi}

\def\resetall{\global\advance\sectno by 1\subsectno=0
   \gdef\firstpart{\number\sectno}\r@s@t}
\def\resetsub{\global\advance\subsectno by 1
   \gdef\firstpart{\number\sectno.\number\subsectno}\r@s@t}
\def\newsection#1\par{\resetall\vskip0pt plus.3\vsize\penalty-250
   \vskip0pt plus-.3\vsize\bigskip\bigskip
   \message{#1}\leftline{\bf#1}\nobreak\bigskip}
\def\subsection#1\par{\ifultr@local\resetsub\fi
   \vskip0pt plus.2\vsize\penalty-250\vskip0pt plus-.2\vsize
   \bigskip\smallskip\message{#1}\leftline{\bf#1}\nobreak\medskip}

\def\t@gsoff#1,{\def\@{#1}\ifx\@\empty\let\next=\relax\else\let\next=\t@gsoff
   \def\@@{p}\ifx\@\@@\else
   \expandafter\gdef\csname#1cite\endcsname##1{\zeigen{##1}}
   \expandafter\gdef\csname#1page\endcsname##1{?}
   \expandafter\gdef\csname#1tag\endcsname##1{\zeigen{##1}}\fi\fi\next}
\def\verbatimtags{\ifx\all\relax\else\expandafter\t@gsoff\all,\fi}
\def\zeigen#1{\hbox{$\langle$}#1\hbox{$\rangle$}}

\def\(#1){\edef\dot@g{\ifmmode\ifinner(\hbox{\noexpand\etag{#1}})
   \else\noexpand\eqno(\hbox{\noexpand\etag{#1}})\fi
   \else(\noexpand\ecite{#1})\fi}\dot@g}

\newif\ifbr@ck
\def\eat#1{}
\def\[#1]{\br@cktrue[\br@cket#1'X]}
\def\br@cket#1'#2X{\def\temp{#2}\ifx\temp\empty\let\next\eat
   \else\let\next\br@cket\fi
   \ifbr@ck\br@ckfalse\br@ck@t#1,X\else\br@cktrue#1\fi\next#2X}
\def\br@ck@t#1,#2X{\def\temp{#2}\ifx\temp\empty\let\neext\eat
   \else\let\neext\br@ck@t\def\temp{,}\fi
   \def\teemp{#1}\ifx\teemp\empty\else\rcite{#1}\fi\temp\neext#2X}
\def\resetbr@cket{\gdef\[##1]{[\rtag{##1}]}}
\def\references{\resetbr@cket\newsection References\par}

\newtoks\symb@ls\newtoks\s@mb@ls\newtoks\p@gelist\n@wcount\ftn@mber
    \ftn@mber=1\newif\ifftn@mbers\ftn@mbersfalse\newif\ifbyp@ge\byp@gefalse
\def\defm@rk{\ifftn@mbers\n@mberm@rk\else\symb@lm@rk\fi}
\def\n@mberm@rk{\xdef\m@rk{{\the\ftn@mber}}%
    \global\advance\ftn@mber by 1 }
\def\rot@te#1{\let\temp=#1\global#1=\expandafter\r@t@te\the\temp,X}
\def\r@t@te#1,#2X{{#2#1}\xdef\m@rk{{#1}}}
\def\b@@st#1{{$^{#1}$}}\def\str@p#1{#1}
\def\symb@lm@rk{\ifbyp@ge\rot@te\p@gelist\ifnum\expandafter\str@p\m@rk=1 
    \s@mb@ls=\symb@ls\fi\write\f@nsout{\number\count0}\fi \rot@te\s@mb@ls}
\def\byp@ge{\byp@getrue\n@wwrite\f@nsin\openin\f@nsin=\jobname.fns 
    \n@wcount\currentp@ge\currentp@ge=0\p@gelist={0}
    \re@dfns\closein\f@nsin\rot@te\p@gelist
    \n@wread\f@nsout\openout\f@nsout=\jobname.fns }
\def\m@kelist#1X#2{{#1,#2}}
\def\re@dfns{\ifeof\f@nsin\let\next=\relax\else\read\f@nsin to \f@nline
    \ifx\f@nline\v@idline\else\let\t@mplist=\p@gelist
    \ifnum\currentp@ge=\f@nline
    \global\p@gelist=\expandafter\m@kelist\the\t@mplistX0
    \else\currentp@ge=\f@nline
    \global\p@gelist=\expandafter\m@kelist\the\t@mplistX1\fi\fi
    \let\next=\re@dfns\fi\next}
\def\symbols#1{\symb@ls={#1}\s@mb@ls=\symb@ls} 
\def\bigsymbol{\textstyle}
\symbols{\bigsymbol\ast,\dagger,\ddagger,\sharp,\flat,\natural,\star}
\def\ftnumbers{\ftn@mberstrue} \def\ftsymbols{\ftn@mbersfalse}
\def\paginal{\byp@ge} \def\resetftnumbers{\ftn@mber=1}
\def\ftnote#1{\defm@rk\expandafter\expandafter\expandafter\footnote
    \expandafter\b@@st\m@rk{#1}}

\long\def\jump#1\endjump{}
\def\ssum{\mathop{\lower .1em\hbox{$\textstyle\Sigma$}}\nolimits}

\def\qed{\nobreak\kern 1em \vrule height .5em width .5em depth 0em}
\def\newneq{\hbox{\rlap{\hbox to 1\wd9{\hss$=$\hss}}\raise .1em 
   \hbox to 1\wd9{\hss$\scriptscriptstyle/$\hss}}}
\def\subsetne{\setbox9 = \hbox{$\subset$}\mathrel{\hbox{\rlap
   {\lower .4em \newneq}\raise .13em \hbox{$\subset$}}}}
\def\supsetne{\setbox9 = \hbox{$\subset$}\mathrel{\hbox{\rlap
   {\lower .4em \newneq}\raise .13em \hbox{$\supset$}}}}

\def\vbar{\mathchoice{\vrule height6.3ptdepth-.5ptwidth.8pt\kern-.8pt}
   {\vrule height6.3ptdepth-.5ptwidth.8pt\kern-.8pt}
   {\vrule height4.1ptdepth-.35ptwidth.6pt\kern-.6pt}
   {\vrule height3.1ptdepth-.25ptwidth.5pt\kern-.5pt}}
\def\f@dge{\mathchoice{}{}{\mkern.5mu}{\mkern.8mu}}
\def\b@c#1#2{{\rm \mkern#2mu\vbar\mkern-#2mu#1}}
\def\b@b#1{{\rm I\mkern-3.5mu #1}}
\def\b@a#1#2{{\rm #1\mkern-#2mu\f@dge #1}}
\def\bb#1{{\count4=`#1 \advance\count4by-64 \ifcase\count4\or\b@a A{11.5}\or
   \b@b B\or\b@c C{5}\or\b@b D\or\b@b E\or\b@b F \or\b@c G{5}\or\b@b H\or
   \b@b I\or\b@c J{3}\or\b@b K\or\b@b L \or\b@b M\or\b@b N\or\b@c O{5} \or
   \b@b P\or\b@c Q{5}\or\b@b R\or\b@a S{8}\or\b@a T{10.5}\or\b@c U{5}\or
   \b@a V{12}\or\b@a W{16.5}\or\b@a X{11}\or\b@a Y{11.7}\or\b@a Z{7.5}\fi}}

\catcode`\X=11 \catcode`\@=12

\expandafter\ifx\csname citeadd.tex\endcsname\relax
\expandafter\gdef\csname citeadd.tex\endcsname{}
\else \message{Hey!  Apparently you were trying to
\string twice.   This does not make sense.} 
\errmessage{Please edit your file (probably \jobname.tex) and remove
any duplicate ``\string\input'' lines} \fi

\def\sciteu{\sciteerror{undefined}}
\def\sciteuphantom{\complainaboutcitation{undefined}}

\def\sciteerror#1#2{{\mathortextbf{\scite{#2}}}\complainaboutcitation{#1}{#2}}
\def\mathortextbf#1{\hbox{\bf #1}}
\def\complainaboutcitation#1#2{%
\vadjust{\line{\llap{---$\!\!>$ }\qquad scite$\{$#2$\}$ #1\hfil}}}

\sectno=-1   
\localtags
\NoBlackBoxes
\define\mr{\medskip\roster}
\define\sn{\smallskip\noindent}
\define\mn{\medskip\noindent}
\define\bn{\bigskip\noindent}
\define\ub{\underbar}
\define\wilog{\text{without loss of generality}}
\define\ermn{\endroster\medskip\noindent}

\define\dbcu{\dsize\bigcup}
\define \nl{\newline}
\documentstyle {amsppt}
\topmatter
\title{On quantification with a finite universe} \endtitle
\author {Saharon Shelah \thanks {\null\newline I would like to thank 
Alice Leonhardt for the beautiful typing. \null\newline
 Publication No. 639 \null\newline
 First Typed - 97/Dec/19 \null\newline
 Latest Revision - 98/Oct/27} \endthanks} \endauthor 
\affil{Institute of Mathematics\\
 The Hebrew University\\
 Jerusalem, Israel
 \medskip
 Rutgers University\\
 Mathematics Department\\
 New Brunswick, NJ  USA} \endaffil
\mn
\abstract   We consider a finite universe ${\Cal U}$ (more exactly - a
family ${\frak U}$ of them), second order quantifiers $Q_K$, where for
each ${\Cal U}$ this means quantifying over a family of $n(K)$-place relations
closed under permuting ${\Cal U}$.  We define some natural orders and shed
some light on the classification problem of those quantifiers. \endabstract
\endtopmatter
\document  

\expandafter\ifx\csname alice2jlem.tex\endcsname\relax
  \expandafter\gdef\csname alice2jlem.tex\endcsname{}
\else \message{Hey!  Apparently you were trying to
\string  twice.   This does not make sense.}
\errmessage{Please edit your file (probably \jobname.tex) and remove
any duplicate ``\string\input'' lines} \fi

\expandafter\ifx\csname bib4plain.tex\endcsname\relax
  \expandafter\gdef\csname bib4plain.tex\endcsname{}
\else \message{Hey!  Apparently you were trying to \string twice.   This does not make sense.}
\errmessage{Please edit your file (probably \jobname.tex) and remove
any duplicate ``\string\input'' lines} \fi

\def\renewcommand{\newcommand}	       
\edef\cite{\the\catcode`@}%
\catcode`@ = 11
\let\@oldatcatcode = \cite
\chardef\@letter = 11
\chardef\@other = 12
%
%
%
%
\def\@innerdef#1#2{\edef#1{\expandafter\noexpand\csname #2\endcsname}}%
%
%
\@innerdef\@innernewcount{newcount}%
\@innerdef\@innernewdimen{newdimen}%
\@innerdef\@innernewif{newif}%
\@innerdef\@innernewwrite{newwrite}%
%
%
%
\def\@gobble#1{}%
%
%
%
\ifx\inputlineno\@undefined
   \let\@linenumber = \empty 
\else
   \def\@linenumber{\the\inputlineno:\space}%
\fi
%
%
%
\def\@futurenonspacelet#1{\def\cs{#1}%
   \afterassignment\@stepone\let\@nexttoken=
}%
\begingroup 
\def\\{\global\let\@stoken= }%
\\ 
\endgroup
\def\@stepone{\expandafter\futurelet\cs\@steptwo}%
\def\@steptwo{\expandafter\ifx\cs\@stoken\let\@@next=\@stepthree
   \else\let\@@next=\@nexttoken\fi \@@next}%
\def\@stepthree{\afterassignment\@stepone\let\@@next= }%
%
%
%
\def\@getoptionalarg#1{%
   \let\@optionaltemp = #1%
   \let\@optionalnext = \relax
   \@futurenonspacelet\@optionalnext\@bracketcheck
}%
%
%
\def\@bracketcheck{%
   \ifx [\@optionalnext
      \expandafter\@@getoptionalarg
   \else
      \let\@optionalarg = \empty
      \expandafter\@optionaltemp
   \fi
}%
\def\@@getoptionalarg[#1]{%
   \def\@optionalarg{#1}%
   \@optionaltemp
}%
%
%
%
\def\@nnil{\@nil}%
\def\@fornoop#1\@@#2#3{}%
\def\@for#1:=#2\do#3{%
   \edef\@fortmp{#2}%
   \ifx\@fortmp\empty \else
      \expandafter\@forloop#2,\@nil,\@nil\@@#1{#3}%
   \fi
}%
\def\@forloop#1,#2,#3\@@#4#5{\def#4{#1}\ifx #4\@nnil \else
       #5\def#4{#2}\ifx #4\@nnil \else#5\@iforloop #3\@@#4{#5}\fi\fi
}%
\def\@iforloop#1,#2\@@#3#4{\def#3{#1}\ifx #3\@nnil
       \let\@nextwhile=\@fornoop \else
      #4\relax\let\@nextwhile=\@iforloop\fi\@nextwhile#2\@@#3{#4}%
}%
%
%
%
\@innernewif\if@fileexists
\def\@testfileexistence{\@getoptionalarg\@finishtestfileexistence}%
\def\@finishtestfileexistence#1{%
   \begingroup
      \def\extension{#1}%
      \immediate\openin0 =
         \ifx\@optionalarg\empty\jobname\else\@optionalarg\fi
         \ifx\extension\empty \else .#1\fi
         \space
      \ifeof 0
         \global\@fileexistsfalse
      \else
         \global\@fileexiststrue
      \fi
      \immediate\closein0
   \endgroup
}%
%
%
%
%
\def\bibliographystyle#1{%
   \@readauxfile
   \@writeaux{\string\bibstyle{#1}}%
}%
\let\bibstyle = \@gobble
%
%
\let\bblfilebasename = \jobname
\def\bibliography#1{%
   \@readauxfile
   \@writeaux{\string\bibdata{#1}}%
   \@testfileexistence[\bblfilebasename]{bbl}%
   \if@fileexists
      \nobreak
      \@readbblfile
   \fi
}%
\let\bibdata = \@gobble
%
%
\def\nocite#1{%
   \@readauxfile
   \@writeaux{\string\citation{#1}}%
}%
\@innernewif\if@notfirstcitation
%
%
\def\cite{\@getoptionalarg\@cite}%
%
%
\def\@cite#1{%
   \let\@citenotetext = \@optionalarg
   \printcitestart
   \nocite{#1}%
   \@notfirstcitationfalse
   \@for \@citation :=#1\do
   {%
      \expandafter\@onecitation\@citation\@@
   }%
   \ifx\empty\@citenotetext\else
      \printcitenote{\@citenotetext}%
   \fi
   \printcitefinish
}%
\def\@onecitation#1\@@{%
   \if@notfirstcitation
      \printbetweencitations
   \fi
   \expandafter \ifx \csname\@citelabel{#1}\endcsname \relax
      \if@citewarning
         \message{\@linenumber Undefined citation `#1'.}%
      \fi
      \expandafter\gdef\csname\@citelabel{#1}\endcsname{%
\strut
\vadjust{\vskip-\dp\strutbox
\vbox to 0pt{\vss\parindent0cm \leftskip=\hsize 
\advance\leftskip3mm
\advance\hsize 4cm\strut\openup-4pt 
\rightskip 0cm plus 1cm minus 0.5cm ?  #1 ?\strut}}
         {\tt
            \escapechar = -1
            \nobreak\hskip0pt
            \expandafter\string\csname#1\endcsname
            \nobreak\hskip0pt
         }%
      }%
   \fi
   \csname\@citelabel{#1}\endcsname
   \@notfirstcitationtrue
}%
%
%
\def\@citelabel#1{b@#1}%
%
%
\def\@citedef#1#2{\expandafter\gdef\csname\@citelabel{#1}\endcsname{#2}}%
%
%
%
\def\@readbblfile{%
   \ifx\@itemnum\@undefined
      \@innernewcount\@itemnum
   \fi
   \begingroup
      \def\begin##1##2{%
         \setbox0 = \hbox{\biblabelcontents{##2}}%
         \biblabelwidth = \wd0
      }%
      \def\end##1{}
      %
      %
      \@itemnum = 0
      \def\bibitem{\@getoptionalarg\@bibitem}%
      \def\@bibitem{%
         \ifx\@optionalarg\empty
            \expandafter\@numberedbibitem
         \else
            \expandafter\@alphabibitem
         \fi
      }%
      \def\@alphabibitem##1{%
         \expandafter \xdef\csname\@citelabel{##1}\endcsname {\@optionalarg}%
         \ifx\biblabelprecontents\@undefined
            \let\biblabelprecontents = \relax
         \fi
         \ifx\biblabelpostcontents\@undefined
            \let\biblabelpostcontents = \hss
         \fi
         \@finishbibitem{##1}%
      }%
      \def\@numberedbibitem##1{%
         \advance\@itemnum by 1
         \expandafter \xdef\csname\@citelabel{##1}\endcsname{\number\@itemnum}%
         \ifx\biblabelprecontents\@undefined
            \let\biblabelprecontents = \hss
         \fi
         \ifx\biblabelpostcontents\@undefined
            \let\biblabelpostcontents = \relax
         \fi
         \@finishbibitem{##1}%
      }%
      \def\@finishbibitem##1{%
         \biblabelprint{\csname\@citelabel{##1}\endcsname}%
         \@writeaux{\string\@citedef{##1}{\csname\@citelabel{##1}\endcsname}}%
         \ignorespaces
      }%
      %
      %
      \let\em = \bblem
      \let\newblock = \bblnewblock
      \let\sc = \bblsc
      \frenchspacing
      \clubpenalty = 4000 \widowpenalty = 4000
      \tolerance = 10000 \hfuzz = .5pt
      \everypar = {\hangindent = \biblabelwidth
                      \advance\hangindent by \biblabelextraspace}%
      \bblrm
      \parskip = 1.5ex plus .5ex minus .5ex
      \biblabelextraspace = .5em
      \bblhook
      \input \bblfilebasename.bbl
   \endgroup
}%
%
%
\@innernewdimen\biblabelwidth
\@innernewdimen\biblabelextraspace
%
%
%
\def\biblabelprint#1{%
   \noindent
   \hbox to \biblabelwidth{%
      \biblabelprecontents
      \biblabelcontents{#1}%
      \biblabelpostcontents
   }%
   \kern\biblabelextraspace
}%
%
%
%
\def\biblabelcontents#1{{\bblrm [#1]}}%
%
%
\def\bblrm{\rm}%
%
%
\def\bblem{\it}%
%
%
\def\bblsc{\ifx\@scfont\@undefined
              \font\@scfont = cmcsc10
           \fi
           \@scfont
}%
%
%
\def\bblnewblock{\hskip .11em plus .33em minus .07em }%
%
%
\let\bblhook = \empty
%
%
%
\def\printcitestart{[}
\def\printcitefinish{]}
\def\printbetweencitations{, }
\def\printcitenote#1{, #1}
%
%
%
\let\citation = \@gobble
%
%
%
\@innernewcount\@numparams
%
%
\def\newcommand#1{%
   \def\@commandname{#1}%
   \@getoptionalarg\@continuenewcommand
}%
%
%
\def\@continuenewcommand{%
   \@numparams = \ifx\@optionalarg\empty 0\else\@optionalarg \fi \relax
   \@newcommand
}%
%
%
\def\@newcommand#1{%
   \def\@startdef{\expandafter\edef\@commandname}%
   \ifnum\@numparams=0
      \let\@paramdef = \empty
   \else
      \ifnum\@numparams>9
         \errmessage{\the\@numparams\space is too many parameters}%
      \else
         \ifnum\@numparams<0
            \errmessage{\the\@numparams\space is too few parameters}%
         \else
            \edef\@paramdef{%
               \ifcase\@numparams
                  \empty  No arguments.
               \or ####1%
               \or ####1####2%
               \or ####1####2####3%
               \or ####1####2####3####4%
               \or ####1####2####3####4####5%
               \or ####1####2####3####4####5####6%
               \or ####1####2####3####4####5####6####7%
               \or ####1####2####3####4####5####6####7####8%
               \or ####1####2####3####4####5####6####7####8####9%
               \fi
            }%
         \fi
      \fi
   \fi
   \expandafter\@startdef\@paramdef{#1}%
}%
%
%
%
%
\def\@readauxfile{%
   \if@auxfiledone \else 
      \global\@auxfiledonetrue
      \@testfileexistence{aux}%
      \if@fileexists
         \begingroup
            \endlinechar = -1
            \catcode`@ = 11
            \input \jobname.aux
         \endgroup
      \else
         \message{\@undefinedmessage}%
         \global\@citewarningfalse
      \fi
      \immediate\openout\@auxfile = \jobname.aux
   \fi
}%
%
%
\newif\if@auxfiledone
\ifx\noauxfile\@undefined \else \@auxfiledonetrue\fi
%
%
%
%
\@innernewwrite\@auxfile
\def\@writeaux#1{\ifx\noauxfile\@undefined \write\@auxfile{#1}\fi}%
%
%
%
\ifx\@undefinedmessage\@undefined
   \def\@undefinedmessage{No .aux file; I won't give you warnings about
                          undefined citations.}%
\fi
%
%
\@innernewif\if@citewarning
\ifx\noauxfile\@undefined \@citewarningtrue\fi
%
%
%
\catcode`@ = \@oldatcatcode


\def\widestnumber#1#2{}

\def\rm{\fam0 \tenrm}

\def\fakesubhead#1\endsubhead{\bigskip\noindent{\bf#1}\par}


%
%
%

%

\font\textrsfs=rsfs10
\font\scriptrsfs=rsfs7
\font\scriptscriptrsfs=rsfs5

\newfam\rsfsfam
\textfont\rsfsfam=\textrsfs
\scriptfont\rsfsfam=\scriptrsfs
\scriptscriptfont\rsfsfam=\scriptscriptrsfs

\edef\oldcatcodeofat{\the\catcode`\@}
\catcode`\@11

\def\Cal@@#1{\noaccents@ \fam \rsfsfam #1}

\catcode`\@\oldcatcodeofat

\newpage

\head {Annotated Content} \endhead  \resetall 
\bn
\S0 $\quad$ Introduction
\mr
\item "{{}}"  [We explain our problem: classifying second order quantifiers
for finite model theory.  We review relevant works, mainly, the work done
on infinite ones.  We then define the basic order relations on such
quantifiers interpretability and expressability.  We also explain why they
are reasonable: as for definable quantifiers, those give the desired
recursiveness result.]
\endroster
\bn
\S1 $\quad$ On some specific quantifiers
\mr
\item "{{}}"  [We define the quantifier we shall use: monadic, partial
one-to-one functions, equivalence relations and linear order.  All have
versions with a cardinality restriction (say cardinality of the domain of
a one-to-one function), which depends on ${\Cal U}$, the universe, only.
For example $Q^{\text{mon}}_{\le \lambda}$ is the quantifier over sets of
cardinality $\le \lambda,Q^{1-1}_{\le \lambda}$ is the quantifier over unary
one-to-one functions with domain of cardinality $\le \lambda$, and 
$Q^{\text{eq}}$ the quantifier over equivalence relations.  We shall 
investigate the natural partial orders on them (by the so-called
interpretability and expressibility).]
\endroster
\bn
\S2 $\quad$ Monadic analyses of $\exists_R$
\mr
\item "{{}}"  [Concentrating first on $\exists_R$, quantifying on the
isomorphic copies of one $n(R)$-place relation $R$, we try to analyze its
``monadic content".  We essentially characterize the maximal cardinality of
a set interpreted by cases of $R$ by a first order formula (actually of low
quantifier depth) as $\lambda_0(R)$ and show that using such a set we can
reduce $R$ to $R_1$ which has domain of cardinality $\lambda_0(R)$.  So up 
to bi-interpretability, $Q_R$ and $\{Q_{R_1},Q_{\le \lambda_0(R)}\}$ are
equivalent.  Now when $\lambda_0(R)$ is too near to the cardinality of the
universe ${\Cal U}$, we have to be more careful but we interpret the (full)
monadic quantifier (with no cardinality restriction).  Lastly, we do the
same for $Q_K$.]
\endroster
\bn
\S3 $\quad$ The one-to-one function analysis
\mr
\item "{{}}"  [We define a cardinal $\lambda_1(R)$ which essentially
characterizes the maximal cardinality of the domain of a one-to-one function
interpretable by cases of $R$.  It is called $\lambda_1(R)$ and we can find
a set $A \subseteq {\Cal U}$ such that the order of magnitude of its 
cardinality is $\lambda_1(R)$ (here - a constant multiple), and show that
$Q_R$ is equivalent by bi-interpretability to $\{Q_{R_1},Q^{\text{eq}}_E\}$
where $E$ is an equivalence relation with not too many equivalence classes
and $R_1$ has domain of cardinality $\sim \lambda_1(R)$.  Of course, $Q_K$
is analyzed similarly. \ub{Now}, unlike the infinite case, up to
bi-expressability $Q^{1-1}_\lambda$ is maximal in the sense that if $R_1$
has domain $\le \lambda^{1/n(R)}$ then it is expressible by $Q^{1-1}_\lambda$.
Hence, under bi-expressibility and up to polynomial order of 
magnitude we have a complete classification.  Of course, on top of 
$Q^{1-1}_{\lambda_1(R)}$ we have the equivalence relation, which is 
understood.]
\endroster
\newpage

\head {\S0 Introduction} \endhead  \resetall 
\bigskip

We investigate and classify to a large extent quantifiers in the following
framework
\mr
\item "{$(*)$}"  for a natural number $n$, 
for a (large) finite set ${\Cal U}$, consider a quantifier $Q_K$ on 
$n$-place relations on ${\Cal U}$, so $K$ is a family of $n$-place relations
on ${\Cal U}$ close under isomorphism (i.e. permutation of ${\Cal U}$).
\ermn
It is natural to restrict ourselves to such families defined by the logic
we have in mind (usually first order), but it seems natural to investigate
two partial orders, interpretability and expressibility defined below, 
which for such
definable classes give the right answer so the use of definability occurs
only in the conclusion.

Earlier this was investigated for infinite ${\Cal U}$, see (below and) in
\cite{Sh:28}, \cite{Bl}, \cite{Sh:171}, but though related, there are some 
differences.  A related work is \cite{BlSh:156} which deals mainly with 
monadic logic on the class of models of a first order theory $T$, so 
its complicatedness measures the complexity of $T$.
We have said on some occasion during this decade that those are
adaptable to finite model theory.  Here we deal with this and shall 
continue in \cite{Sh:F334}.  

In \cite{Sh:28} we gave a complete classification of the class of second order
quantifiers: those which are first-order definable (see below an exact
definition).  We find that for infinite models up to a very strong notion of
equivalence, bi-interpretability, there are only four such quantifiers:
first order, monadic, one-to-one partial functions, and second-order.  See
Baldwin \cite{Bl}.

Now \S1-\S3 of the present work are parallel to \S1, \S2, \S3 of 
\cite{Sh:171}, so below we describe the latter and then explain what we shall
do here.  
In \cite{Sh:171} our aim was to see what occurs if we remove the restriction
that the quantifier is first-order definable.  As we do not want to replace
this by a specific ${\Cal L}$-definable (${\Cal L}$-some logic) we restricted
ourselves in \cite{Sh:171} to a fixed infinite universe ${\Cal U}$.  If 
we then want to restrict
ourselves to ${\Cal L}$-definable quantifiers, we are able to remove the
restriction to a fixed universe ${\Cal U}$.

The strategy in \cite{Sh:171} is to squeeze the quantifier $Q_R$ (similarly
for $Q_K$) between some well understood quantifiers to get, eventually, 
equality.  Unfortunately, for interpretability we get a lower bound and 
an upper bound which are close but 
not necessarily equal; i.e. both of the form $Q_{\bold E}$, where $\bold E$
is a set of equivalence relations and they are quite close (see below).  
More specifically we use cases of
$Q^{\text{eq}}_{\lambda,\mu}$ (i.e. on equivalence relations with $\lambda$
classes each of cardinality $\le \mu$).  Carrying out the strategy we first 
``find" the monadic content of, say, $Q_R$, by interpreting in it 
$Q^{\text{mon}}_{\lambda_0(R)}$ which is quantifying on sets of cardinality
$\le \lambda_0(R)$ and $\lambda_0(R)$ is maximal (and reduce the problem to
``the remainder", that is a relation $R_1$ with Dom$(R_1)$ of cardinality
$\le \lambda_0(R)$ and $Q_{R_1} \le_{\text{int}} Q_R$).  Next interpret
$Q^{1-1}_{\lambda_1(R)}$ which is quantifying on partial one-to-one
functions of cardinality $\le \lambda_1(R)$.  Now we succeed to squeeze $Q_R$,
for ``the remainder" between $Q_{\lambda,\lambda}$ and 
$Q_{\mu,\mu},\lambda \le \mu \le \text{ Min}\{2^\lambda,|{\Cal U}|\}$
but in general cannot show this with $\lambda = \mu$.  Clearly if 
$|{\Cal U}|$ is $\aleph_0$, this does not occur and we can get a complete
picture (see below \scite{0.1}).  Also by ``expressibility" 
(a stronger equivalence relation but O.K. for the application to logic) 
if $V=L$, then the gap does not occur, but in some generic extensions it does.
\bn
So by \cite{Sh:171} we can e.g. conclude
\proclaim{\stag{0.B} Theorem}  Assume $K$ is a family of $n$-place relations 
over ${\Cal U}$ where $|{\Cal U}| = \aleph_0$.  \ub{Then} $Q_K$ is 
bi-interpretable (see below) with $Q_{\bold E}$ for some family $\bold E$ 
of equivalence relations.
\endproclaim
\bn
We can make this more specific.

The present situation is more complicated.  For example, the finite
cardinalities allow a family of monadic quantifiers: for the case
$|{\Cal U}| = n$ we have $Q_{\text{\rm ln n\/}},Q_{\text{\rm ln ln n\/}}$, 
etc.  However, modulo these
cardinality restrictions we are able to get a picture analogous to the 
original case.  Also in the fine analysis we do not get an equivalence
relation $E$ on ${\Cal U}$ such that $Q_R,Q_E$ are bi-interpretable or even
just bi-expressible, but just ``squeeze" $Q_R$ between two such quantifiers,
which are quite closed (i.e. size of one bounded by polynomial in the size
of another).  That is (concentrating on the case ${\Cal U}$ is fixed (and
finite)): assume $R$ is an $n$-place relation on ${\Cal U}$ then we can
\ub{uniformly} attach it to a cardinal $\lambda_1(R)$, and an equivalence
relation $E$ such that:
\mr
\item "{$(\alpha)$}"  $Q^{\text{eq}}_E,Q^{1-1}_{\lambda_1(R)}$ are
interpretable in $Q_R$ (quantifiers over equivalence relations isomorphic
to $E$ and partial $1-1$-functions of cardinality $\le \lambda_1(R))$
\sn
\item "{$(\beta)$}"  if $\lambda = \lambda_1(R)^{n(R)} \le |{\Cal U}|$
\ub{then} $Q_R$ is expressible by $(Q^{\text{eq}}_E,Q^{1-1}_\lambda)$
\sn
\item "{$(\gamma)$}"  if $\lambda_1(R)^{n(R)} > |{\Cal U}|$, \ub{then} any
binary relation on a set $A \subseteq U$ with cardinality 
$|A|^{\frac{1}{2n(R)}}$ is interpretable in $Q_R$.
\ermn
The uniformly means that the formulas involved in interpretability or
expressibility does not depend on $R$ and ${\Cal U}$ but on $n$, in fact we
can give explicit bounds on their size from $n$.

Note that we abuse notation using $R$ as a relation and predicate; of course,
the formulas have an $n$-place predicate to stand for copies of $R$ (see
below). 

Note we actually deal also with quantifying on appropriate families of $R$'s
of fix arity (e.g. those satisfying some sentence).  Note that we cannot get
much better results by counting.
\bn
\centerline{$* \qquad * \qquad *$}
\bn
We thank C. Steinhorn, J. Tyszkiewicz and J. Baldwin for helpful discussions
on preliminary versions in MSRI 10/89, Dimacs 95/96 and Rutgers Fall 1997, 
respectively.  Much more is due to Baldwin, Fall 1998, for helping to greatly
improve the presentation.  
\bn
Let us now make some conventions and definitions.
\demo{\stag{0.1} Convention}  1) Informally ${\Cal U}$ will be a fixed 
finite universe
(usually large compared to $n$) but, if not said otherwise, we are proving
things uniformly.  So more exactly, ${\Cal U}$ varies on ${\frak U}$, a 
family of such sets.  You may choose ${\frak U} = \{(0,n):n \text{ a natural
number}\}$.  \nl
2)  Informally, $K$ will denote a family of 
$n$-place relations over ${\Cal U}$, (for a
natural number $n=n(K))$, closed under isomorphism, i.e. if $R_1,R_2$ are
$n$-place relations on ${\Cal U}$ and $({\Cal U},R_1) \cong 
({\Cal U},R_2)$ \ub{then} $R_1 \in K$ iff $R_2 \in K$.  So \ub{formally} 
$K$ is a function with domain ${\frak U}$ and $K[{\Cal U}]$ is as above; but
$n(K) = n(K[{\Cal U}])$ for each ${\Cal U} \in {\frak U}$. Also below without
saying in e.g. Definition \scite{0.4} the formula $\varphi$ is the same for 
all ${\Cal U} \in {\frak U}$.  \nl
3) Let $\bar K$ denote a finite sequence of such $K$'s, that is

$$
\bar K = \langle K_\ell:\ell < \ell g(\bar K) \rangle, \text{ so } 
\bar K^j_i = \langle K^j_{i,\ell}:\ell < \ell g(\bar K^j_i) \rangle.
$$
\mn
4) Let $R$ denote a relation, its domain is
Dom$(R) = \cup\{\bar a: \models R(\bar a)\}, n=n(R)$ if $R$ is an $n$-place
relation (or predicate; we shall not always strictly distinguish).  
Usually $R$ is on ${\Cal U}$ which is clear from the context.  
\ub{Formally}, $R$ is a function with domain ${\frak U}$ and 
$R[{\Cal U}]$ is an $n(R)$-place relation on ${\Cal U}$).
\enddemo
\bigskip

\definition{\stag{0.2} Definition}  For any $K,\exists_K$ (or $Q_K$) denotes 
a second order quantifier, intended to vary on members of $K$.  More exactly,
$L(\exists_{K_1},\dotsc,\exists_{K_m})$ is defined like first order
logic but we have for each $\ell = 1,m$ (infinitely many) variables $R$ which
serve as $n(K_\ell)$-place predicates, and we can form $(\exists_{K_i}R)
\varphi$ for a formula $\varphi$ (when $R$ is $n(K_i)$-place).  
Defining satisfaction, we look only 
at models with universe ${\Cal U}$, and $\models (\exists_{K_\ell}R)\varphi
(R,\cdots)$ iff for some $R^0 \in K_\ell[{\Cal U}]$ we 
have $\varphi(R^0,\cdots)$.

We may display the predicates (or relations) appearing in $\varphi$, i.e.
$\varphi(x,y,\bar R)$.  Of course, we may write $K$ not $K[{\Cal U}]$, etc.,
abusing notation.
\enddefinition
\bigskip

\remark{Remark}  Note that quantifiers depending on parameters are not
allowed, e.g. automorphisms; on such quantifiers see \cite{Sh:e}.
\endremark
\bigskip

\definition{\stag{0.3} Definition}  We say that $K$ (or $Q_K$) is 
${\Cal L}$-definable (where ${\Cal L}$ is a logic) if there is a 
formula $\varphi(R) \in {\Cal L}$, in the vocabulary $\{R\}$ and is 
appropriate, i.e. an $n(K)$-place
predicate, such that for any $n$-place relation $R$ on ${\Cal U}$

$$
({\Cal U},R) \models \varphi(R) \text{ iff } R \in K.
$$
\enddefinition
\bigskip

\definition{\stag{0.4} Definition}  1) We say that 
$\exists_{K_1} \le_{\text{int}} \exists_{K_2}$ (in other words 
$\exists_{K_1}$ is interpretable in $\exists_{K_2}$) \ub{if} for some 
first-order formula $\varphi(\bar x,\bar S) = \varphi(x_0,\dotsc,x_{n(K_1)-1},
S_0,\dotsc,S_{m-1})$, (each $S_\ell$ is an
$n(K_2)$-place predicate) the following holds:
\mr
\item "{$(*)$}"  for every ${\Cal U} \in {\frak U}$ and $R_1 \in 
K_1[{\Cal U}]$ there are $S_0,\dotsc,S_{m-1} \in K_2[{\Cal U}]$ such that
$({\Cal U},S_0,\dotsc,S_{m-1}) \models (\forall \bar x)[R_1(\bar x) \equiv
\varphi(\bar x,S_0,\dotsc,S_{m-1})]$ 
\ermn
(so in $(*)$, $\varphi$ does not depend on ${\Cal U}$). \nl
2) We say $k$-interpretable \ub{if} we demand $m \le k$, and then write
$\le_{k\text{-int}}$. \nl
3) We can define $\exists_{K_1} \le^{\Cal L}_{\text{int}} \exists_{K_2}$ or
$\exists_{K_1} \le_{\text{int}} \exists_{K_2}$ mod ${\Cal L}$ similarly, by 
letting $\varphi \in {\Cal L}$.  Similarly for
$\le^{\Cal L}_{k\text{-int}}$.  Instead we may say modulo ${\Cal L}$.
\enddefinition
\bn
We define a weaker relative of interpretability; we say $\exists_{K_1}$ is
expressible by $\exists_{K_2}$ if in the notion of interpretable we take
the formula $\varphi$ to be in the logic $L(\exists_{K_2})$.  This is then
a special but very important case of \scite{0.4}(3).
\definition{\stag{0.5} Definition}  1) We say that $\exists_{K_1}
\le_{\text{exp}} \exists_{K_2}$ (in other words $\exists_{K_1}$ is expressible
by $\exists_{K_2}$) \ub{if} there is a formula $\varphi(\bar x,S_0,\dotsc,
S_{m-1})$ in the logic $L(\exists_{K_2})$ such that:
\mr
\item "{$(*)$}"  for every ${\Cal U} \in {\frak U}$ and $R_1 \in 
K_1[{\Cal U}]$, there are $S_0,\dotsc,S_{m-1} \in K_2[{\Cal U}]$ such that
$({\Cal U},S_0,\dotsc,S_{m-1}) \models (\forall \bar x)[R_1(\bar x) \equiv
\varphi(\bar x,S_0,\dotsc,S_{m-1})]$.
\ermn
2) We say that $\exists_{K_1} \le_{\text{inex}} \exists_{K_2}$ 
(in other words $\exists_{K_2}$ is invariantly expressible by 
$\exists_{K_2}$) \ub{if} there is a formula $\varphi
(\bar x,S_0,\dotsc,S_{m-1})$ in the logic $L(\exists_{K_2})$ such that:
\mr
\item "{$(*)$}"  for every ${\Cal U} \in {\frak U}$ and $R_1 \in K_1
[{\Cal U}]$, there are $S_0,\dotsc,S_{m-1} \in K_2[{\Cal U}]$ such that 
for every $K_3$ which extends $K_2$, letting $\varphi'$ be 
$\varphi$ when we replace $\exists_{K_2}$ by $\exists_{K_3}$ we have:
\endroster

$$
({\Cal U},S_0,\dotsc,S_{m-1}) \models (\forall \bar x)[R_1(\bar x) \equiv
\varphi'(\bar x,S_0,\dotsc,S_{m-1})].
$$
\mn
3) We define $k$-expressible, $\le_{k\text{-exp}}$, invariantly
$k$-expressible and $\le_{k\text{-inex}}$ and may add ${\Cal L}$ as a
superscript parallel to \scite{0.4}(2).
\enddefinition
\bigskip

\definition{\stag{0.6} Definition}  1) We say that $\exists_{K_1} 
\equiv_{\text{int}} \exists_{K_2}$ 
(in other words $\exists_{K_1},\exists_{K_2}$ are bi-interpretable) if 
$\exists_{K_1} \le_{\text{int}} \exists_{K_2}$ and $\exists_{K_2} 
\le_{\text{int}} \exists_{K_1}$. \nl
2) We say $\exists_{K_2} \equiv_{\text{exp}} \exists_{K_2}$ 
(in other words $\exists_{K_1},\exists_{K_2}$ are bi-expressible) if 
$\exists_{K_1} \le_{\text{exp}} \exists_{K_2}$ and $\exists_{K_2} 
\le_{\text{exp}} \exists_{K_1}$.  Similarly for $\equiv_{\text{inex}}:
\exists_{K_1} \equiv_{\text{inex}} \exists_{K_2}$ (in other words
$\exists_{K_2},\exists_{K_1}$ are invariantly bi-expressible) if
$\exists_{K_1} \le_{\text{inex}} \exists_{K_2}$ and $\exists_{K_2}
\le_{\text{inex}} \exists_{K_1}$. \nl
3) We can define $\exists_{K_1} \le_{\text{int}} \{\exists_{K_0},\dotsc,
\exists_{K_{k-1}}\}$ as in Definition \scite{0.4} but $S_0,\dotsc, \in
\dbcu^{k}_{i=1} K_i[{\Cal U}]$, we let 
$\exists_{\bar K}$ stand for $\{\exists_{K_0},
\dotsc,\exists_{K_{k-1}}\}$ where $K = \langle K_0,\dotsc,K_{k-1} \rangle$;
we define $\exists_{\bar K^1} \le_{\text{int}} \exists_{\bar K^2}$ if
$\exists_{K^1_\ell} \le_{\text{int}} \exists_{\bar K^2}$ for each $\ell$;
we also define expressible, invariantly expressible, bi-interpretable and
(invariantly) bi-expressible similarly. \nl
4) Let $\exists_{K_1} \equiv_{1\text{-int}} \exists_{K_2}$ (in other words
$\exists_{K_1},\exists_{K_2}$ are 1-bi-interpretable) if $\exists_{K_1}
\le_{1\text{-int}} \exists_{K_2}$ and $\exists_{K_2} \le_{1\text{-int}}
\exists_{K_1}$; recall $\le_{1-\text{int}}$ is defined in \scite{0.4}(2) for
$k=1$.  Similarly $\exists_{K_1} \equiv_{1\text{-exp}} \exists_{K_2}$
and $\exists_{K_1} \equiv_{1\text{-inex}} \exists_{K_2}$. \nl
5) In all those notions we add ``modulo $\bar K$" if parameters from
$\cup\{K_\ell:\ell < \ell g(\bar K)\}$ are allowed.  We can combine this with
\scite{0.4}(3) so have modulo $(\bar K,{\Cal L})$.
\enddefinition
\bigskip

\demo{\stag{0.7} Notation}  1) If $R_\ell$ is an $n_\ell$-place relation
for $\ell < n$ then we let
$\dsize \sum_{\ell=0}^{n-1} R_\ell = \{\bar a_0 \char 94 \cdots \char 94
\bar a_{n-1}:\bar a_\ell \in R_\ell\}$; more formally
$(\dsize \sum^{n-1}_{\ell = 0} R_\ell)({\Cal U}) = \dsize \sum^{n-1}
_{\ell=0} R_\ell[{\Cal U}]$. \nl
2) Let $\dsize \sum^{n-1}_{\ell=0} K_\ell = \{\dsize \sum^{n-1}_{\ell=0}
R_\ell:R_\ell \in K_\ell \text{ for } \ell < n\}$. \nl
3) $\exists_R$ stands for $\exists_K$ where $K = \{R_1:({\Cal U},R_1) \cong
({\Cal U},R)\}$ and so formally if $R = \langle R[{\Cal U}]:{\Cal U} \in
{\frak U} \rangle$ then $K_R$ is defined by $K_R[{\Cal U}] = \{R_1:
({\Cal U},R_1) \cong ({\Cal U},R)\}$.
\enddemo
\bigskip

\proclaim{\stag{0.8} Lemma}  1) $\le_{\text{int}},\le_{\text{inex}}$ and
$\le_{\text{exp}}$ as well as $\le_{1\text{-int}},\le_{1\text{-inex}}$ and
$\le_{1\text{-exp}}$ are partial quasi orders.  Hence $\equiv_{\text{int}},
\equiv_{\text{inex}},\equiv_{\text{exp}}$ are equivalence relations
as well as $\equiv_{1\text{-int}},\equiv_{1\text{-inex}}$ and
$\equiv_{1\text{-exp}}$. \nl
2) $\exists_{\bar K_1} \le_{\text{int}} \exists_{\bar K_2}$ implies
$\exists_{\bar K_1} \le_{\text{inex}} \exists_{\bar K_2}$ which implies
$\exists_{\bar K_1} \le_{\text{exp}} \exists_{\bar K_2}$. 
Similarly for the ``1-" versions.  Also each ``1-" version implies the one
without.  \nl
3) $\exists_{\bar K}$ and $\exists_K$ are bi-interpretable if $K = \dsize
\sum_i K_i$ or $K = \dbcu_i K_i$ (where 
$n(K_i)$ constant in the second case). \nl
4) In all those cases we can do everything modulo $\bar K_0$ or modulo
${\Cal L}$ (if ${\Cal L}$ is a reasonable logic closed by first order
operations) or modulo $(\bar K_0,{\Cal L})$.
\endproclaim
\bigskip

\demo{Proof}  Straight.
\enddemo
\bigskip

\proclaim{\stag{0.9} Lemma}  1) If $\bar K_1,\bar K_2$ are ${\Cal L}$-definable
(i.e. each $K_{\ell,i}$ is, see Definition \scite{0.3}) and 
$\exists_{\bar K_1} \le_{\text{exp}}
\exists_{\bar K_2}$ \ub{then} we can recursively attach to every formula in
${\Cal L}(\exists_{\bar K_1})$ an equivalent formula in ${\Cal L}
(\exists_{\bar K_2})$. \nl
2) If $\bar K_1,\bar K_2$ are ${\Cal L}$-definable, $\exists_{\bar K_1}
\le_{\text{exp}} \exists_{\bar K_2}$ \ub{then} the set of valid ${\Cal L}
(\exists_{\bar K_1})$-sentences that is ${\Cal L}(\exists_{K_{1,0}},\dotsc,
\exists_{K_{1,\ell g(\bar K)-1}})$-sentences, recursive in the set of valid
${\Cal L}(\exists_{\bar K_2})$-sentences.
\endproclaim
\bigskip

\demo{Proof}  Easy.
\enddemo
\bigskip

\remark{Remark}   1) The need of ``${\Cal L}$-definable" is clearly necessary.
Though at first glance the conclusions of \scite{0.9} may seem the natural
definition of interpretable, I think reflection will lead us to see it isn't.
\nl
2)  Note that naturally we use \scite{0.9} with \scite{0.8}. \nl
3)  Note that, of course, in \scite{0.9}, it is understood that
the formulas from ${\Cal L}$ are the same for all ${\Cal U} \in {\frak U}$\,.
\endremark
\newpage

\head {\S1 On some specific quantifiers} \endhead  \resetall 
\bigskip

\definition{\stag{1.1} Definition}  0) $K^{\text{tr}} = \{A \subseteq
{\Cal U}:|A|=1\}$, and we can write $\exists$ for $\exists_{K^{\text{tr}}}$;
here tr stands for trivial.
\nl
1) Let $K^{\text{mon}}_\lambda = \{A \subseteq {\Cal U}:|A| = \lambda\}$ 
for a number $\lambda \le |{\Cal U}|/2$; here mon stands for monadic. \nl
2) But we write $Q^{\text{mon}}_\lambda$ for 
$\exists_{K^{\text{mon}}_\lambda}$, and similarly for the other quantifiers
defined below. \nl
3) $K^{1-1}_\lambda = \{f:f$ is a partial one-to-one function,
$|\text{Dom}(f)| = \lambda\}$ when $\lambda \le |{\Cal U}|/2$. \nl
4) $K^{\text{eq}}_{\lambda,\mu} = \{E:E$ is an equivalence relation on some
$A \subseteq {\Cal U}$, with $\lambda$ equivalence classes, each of power
$\mu\}$. \nl
5) In 4) we can replace ``$\mu$" by ``$< \mu$" if each equivalence class has
$< \mu$ elements.  Similarly replacing $\lambda$ by ``$< \lambda$".
Similarly $\le \lambda,\le \mu$. \nl
6) $K^{\text{mon}}_{< \lambda} = \dbcu_{\mu < \lambda} K^{\text{mon}}_\mu$
and $K^{1-1}_{< \lambda} = \dbcu_{\mu < \lambda} K^{1-1}_\mu$.  Similarly with
$\le \lambda$; here the ``less than half" is not so important. \nl
7) $K^{\text{mon}} = K^{\text{mon}}_* = \{A:A \subseteq {\Cal U}\}$ and
$K^{1-1} = K^{1-1}_* = \{f:f$ is a partial one-to-one function$\}$. \nl
8) $K^{\text{eq}}_{\lambda,*} = \{E:E$ is an equivalence relation on some
$A \subseteq {\Cal U}$ with $\lambda$-equivalence classes$\}$ and
$K^{\text{eq}}_{*,< \mu} = \{E:E$ is an equivalence relation on some $A
\subseteq {\Cal U}$ each equivalence class $< \mu\}$ and $K^{\text{eq}} =
K^{\text{eq}}_{\lambda,*} = \{E:E$ is an equivalence relation on some $A
\subseteq {\Cal U}\}$ and lastly 
$K^{\text{eq}}_{\le \lambda} = \{E:E$ an equivalence relation on 
$A \subseteq {\Cal U},|A| \le \lambda\}$.
\enddefinition
\bigskip

\remark{\stag{1.1A} Remark}  More formally $\lambda,\mu$, etc., are functions
from ${\frak U}$ to $\Bbb N$ which satisfy conditions such as $\lambda
({\Cal U}) < \frac{|{\Cal U}|}{2}$.  
We write $\frac{|{\Cal U}|}{2}$ as shorthand for
$[\frac{{\Cal U}}{2}]$.

Claims \scite{1.2} through \scite{1.6} are established by
similar arguments.  To illustrate the technique we prove \scite{1.2}(4).  
If $\varphi(x,S_0,S_1)$ denotes $``x \in S_0 \vee x \in S_1"$ then
as $S_1,S_2$ range over subsets of ${\Cal U}$ with $|{\Cal U}| < \lambda <
\frac{|{\Cal U}|}{2}$ clearly all sets of cardinality $\kappa,
\frac{|{\Cal U}|}{4}
\le \kappa \le \frac{|{\Cal U}|}{2}$ are represented; all sets of cardinality
$< \frac{|{\Cal U}|}{4}$ are represented by $x \in S_0 \and x \in S_1$.
(Note this depends on $|S_i| \le \frac{|{\Cal U}|}{2}$).  Finally sets with
cardinality between $\frac{|{\Cal U}|}{2}$ and $|{\Cal U}|$ are represented
by taking compliments.

The choice $\frac{|{\Cal U}|}{2}$ and $\frac{|{\Cal U}|}{4}$ is arbitrary.
But if $\frac{|{\Cal U}|}{k}$ for larger $k$ were choice, the union of two
sets would have to be replaced by a union of more sets.  A lower bound of the
form $\frac{|{\Cal U}|}{k}$ permits the uniform choice of the formula 
$\varphi$.
\endremark
\bigskip

\proclaim{\stag{1.2} Claim}  Let 
$\lambda \le \chi < \frac{|{\Cal U}|}{2}$.  Then, uniformly (the choice of
the interpreting formula $\varphi$ does not depend on ${\Cal U}$) we have: 
\nl
0) $Q^{\text{mon}}_\chi$ is $\exists_R$ for some $R$. \nl
1) $Q^{\text{mon}}_\lambda \equiv_{\text{int}} Q^{\text{mon}}_{\le \lambda}$
and $Q^{\text{mon}}_{< \lambda} \le_{\text{int}} Q^{\text{mon}}_{< \chi}$. \nl
2) $Q^{\text{mon}}_{< \mu} \equiv_{\text{int}} Q^{\text{eq}}_{1,< \mu}$. \nl
3) $Q^{\text{mon}}_\lambda \equiv_{\text{int}} 
Q^{\text{mon}}_{\lambda + \lambda}$ if $\lambda({\Cal U}) \le 
\frac{|{\Cal U}|}{4}$. \nl
4) If $|{\Cal U}|/2 \ge \lambda \ge |{\Cal U}|/4$, \ub{then}
$Q^{\text{mon}}_\lambda \equiv_{\text{int}} Q^{\text{mon}}$.  \nl
5) More generally, for any constants $a$ and $b$, if $\frac{|{\Cal U}|}{a}
\ge \lambda \ge \frac{|{\Cal U}|}{b}$, then $Q^{\text{mon}}_\lambda 
\equiv_{\text{int}} Q^{\text{mon}}$.
\endproclaim
\bigskip

\demo{Proof}  Straightforward.  For 0) recall Notation \scite{0.7}(3).
For 2) recall Definition \scite{1.1}(4).
\enddemo
\bigskip

\proclaim{\stag{1.3} Claim}  Let $\lambda \le \chi$ be as in \scite{1.2}. \nl
0) $Q^{1-1}_\lambda$ is $\exists_R$ for some $R$.
\nl
1) $\lambda \le |{\Cal U}|/2 \Rightarrow Q^{1-1}_\lambda 
\equiv_{\text{int}} Q^{1-1}_{\le \lambda}$; \ub{and} $\chi \le |{\Cal U}|/2
\Rightarrow Q^{1-1}_{< \lambda} \le_{\text{int}} Q^{1-1}_{< \chi}$. \nl
2) $Q^{\text{mon}}_{< \lambda} \le_{\text{int}} Q^{1-1}_{< \lambda}$. \nl
3) If $\lambda \ge |{\Cal U}|/4$, \ub{then} $Q^{1-1}_\lambda
\equiv_{\text{int}} Q^{1-1}_*$. \nl
4) $Q^{1-1}_{\lambda + \lambda} \equiv_{\text{int}} Q^{1-1}_\lambda$ if
$\lambda \le |{\Cal U}|/4$. \nl
5) If $R$ is a graph of a partial one-to-one function on ${\Cal U},\lambda =
\text{ Min}\{|\text{Dom}(R),R\}$ then $Q_R \equiv_{\text{int}} 
Q^{1-1}_\lambda$.
\endproclaim
\bigskip

\demo{Proof}  Straightforward.
\enddemo
\bigskip

\proclaim{\stag{1.4} Claim}  Let $\lambda \le \chi$ and $\mu \le \kappa$
(as in \scite{1.2}). \nl
0) $Q^{\text{eq}}_{\lambda,\mu}$ is $\exists_R$ for some $R$. \nl
1) If $\chi \le |{\Cal U}|/2,\kappa \le |{\Cal U}|/2$, \ub{then}
$Q^{\text{eq}}_{\lambda,\mu} \le_{\text{int}} Q^{\text{eq}}_{\chi,\kappa}
\equiv_{\text{int}} Q^{\text{eq}}_{\le \chi,\le \kappa}$.  \nl
1A)  $Q^{\text{eq}}_{< \lambda,< \mu} \le_{\text{int}} 
Q^{\text{eq}}_{< \chi,< \kappa}$. \nl
2) If $|{\Cal U}|/4 \le \lambda \le |{\Cal U}|/2$ then $Q^{\text{eq}}
_{\le \lambda} \equiv_{\text{int}} Q^{\text{eq}}_{*,*}$. \nl
3) For equivalence relations $E_1,E_2$ on ${\Cal U}$, natural sufficient
condition for interpretability works.  Similarly for families of equivalence
relations. \nl
4) $\exists_K \le_{\text{int}} Q^{\text{eq}}_{\le \lambda}$ if $(\forall R \in
K)[|\text{Dom}(R)|^{n(R)} \le (\lambda -1)^2]$ and $\lambda \le
|{\Cal U}|/2$.
\endproclaim
\bigskip

\demo{Proof}  Left to the reader.
\enddemo
\bigskip

\definition{\stag{1.4A} Definition}  1) $Q^{\text{ord}}_\lambda = \{R:R
\text{ a linear order of a subset } A \text{ of } {\Cal U}$ of cardinality
\footnote{for the infinite case we demand otp$(A,R) = \lambda$}
$\lambda\}$. \nl
2) $Q^{\text{ord}}_{< \lambda} = \dbcu_{\mu < \lambda} Q^{\text{ord}}_\mu$.
\enddefinition
\bigskip

\proclaim{\stag{1.5} Claim}  0) $Q^{\text{ord}}_\lambda$ has the form $Q_R$.
\nl
1) $Q^{\text{mon}}_{\le \lambda} \le_{\text{int}} Q^{\text{ord}}_\lambda$. \nl
2) If $\mu \times \kappa \le \lambda$ then $Q^{\text{eq}}_{\mu,\kappa}
\le_{\text{int}} Q^{\text{ord}}_\lambda$ and $Q^{\text{eq}}_{\le \lambda}
\le_{\text{int}} Q^{\text{ord}}_\lambda$. \nl
3) $\mu < \kappa \le \lambda \Rightarrow Q^{\text{eq}}_{\mu,\kappa}
\le_{1\text{-int}} Q^{\text{ord}}_\lambda$ mod $Q^{\text{ord}}_{\le \lambda}$
and $Q^{\text{eq}}_{\le \lambda} \le_{1\text{-int}} Q^{\text{ord}}_\lambda$
mod $Q^{\text{ord}}_{\le \lambda}$. \nl
4) $Q^{\text{ord}}_\lambda \le_{1\text{-int}} Q^{1-1}_\lambda$ mod
${\Cal L}(Q^{1-1}_{\le \mu})$ if $\lambda \le \mu$. \nl
5) $Q^{\text{ord}}_\lambda \le_{\text{int}} Q^{\text{eq}}_{\lambda,\lambda}$,
in fact, one $E_0 \in Q^{\text{eq}}_{\lambda,\lambda}$, one $E_1 \in
Q^{\text{eq}}_{(\lambda^2),2}$ and one 
$P \in Q^{\text{mon}}_{(\lambda^2)}$ suffice.
\endproclaim
\bigskip

\demo{Proof}  Straight.
\enddemo
\bigskip

\proclaim{\stag{1.6} Claim}  1) $Q_{K_1} \le_{\text{int}} Q_{K_2}$ mod
$Q^{1-1}_*$ is equivalent to $Q_{K_1} \le_{1\text{-int}} Q_{K_2}$ mod
$Q^{1-1}_*$. \nl
2) Similarly for $\le_{\text{inex}},\le_{\text{exp}}$.
\endproclaim
\bigskip

\definition{\stag{1.7} Definition}  For any equivalence relation $E$ on a 
set Dom$(E) \subseteq {\Cal U}$ we define \nl
1) nu$_{\ge k}(E)$ is the number of equivalence classes of $E$ with $\ge k$
members. \nl
2) uq$_k(E) = \text{ Max}\{|B|:B \subseteq {\Cal U} \text{ and there are }
E_0,\dotsc,E_{k-1} \in Q_E[{\Cal U}]$ such that: $b \ne c \in B \Rightarrow
(\exists \ell < k)[(b E_\ell b \equiv \neg c E_\ell c) \vee (b E_\ell b \and 
c E_\ell c \and \neg b E_\ell c)]\}$.
\nl
3) For $x \in {\Cal U} \backslash \text{ Dom}(E)$ let $x/E$ be ${\Cal U}
\backslash \text{ Dom}(E)$.
\enddefinition
\bigskip

\proclaim{\stag{1.8} Claim} 1) $Q^{1-1}_{\text{nu}_{\ge 2}(E)} 
\le_{\text{int}} Q_E$. \nl
2) $Q^{1-1}_{\text{uq}_k(E)} \le_{\text{int}} Q_E$.
\endproclaim
\bigskip

\demo{Proof}  Let ${\Cal U} \in {\frak U}$. \nl
1) We can find a sequence $\langle a_i:i < 2 \text{ nu}_{\ge 2}(E) \rangle$
with no repetitions, $a_i \in {\Cal U}$ such that for $i < j$ we have
$a_i E a_j \Leftrightarrow j = i+1 \and$ ``$i$ is even".  
Let $P_0 = \{a_{2i}:i\},
P_1 = \{a_{2i+1}:i\}$.  So $P_0(x) \and P_1(y) \and x Ey$ defines a partial
one to one function with domain of cardinality $2$ nu$_{\ge 2}(E)$.  We
finish as we can interpret $Q^{\text{mon}}_{|P_0|}$ (or see \scite{2.1A}(2)). 
\nl
2) Easy, too.  \hfill$\square_{\scite{1.8}}$
\enddemo
\bigskip

\definition{\stag{1.11} Definition}  $Q^{n\text{-ary}}_\mu$ is 
quantifying on $n$-place relation with domain of cardinality $\le 
\mu({\Cal U})$.
\enddefinition
\bigskip

\proclaim{\stag{1.12} Claim}  1) For $n$, letting $\bar x = \langle x_0,
\dotsc,x_{n-1}\rangle$ there is a formula $\varphi(\bar x,F_0,\dotsc,
F_{n-1})$ in monadic logic ($F_\ell$ unary function symbol), such that:
\mr
\item "{$(*)$}"  for ${\Cal U} \in {\frak U},A \subseteq {\Cal U}$, an
$n$-place relation $R$ on $A$ we can find a model $M = ({\Cal U},F^M_0,
\dotsc,F^M_{n-1})$ and partial one-to-one functions $F^M_0,\dotsc,F^M_{n-1}$
from ${\Cal U}$ to ${\frak U}$ such that $\varphi(\bar x;F^M_0,\dotsc,
F^M_{n-1})$ define $R$ in $M$, where the monadic quantifier is being 
interpreted as $Q^{\text{mon}}_{< \lambda},\lambda \ge |R|$ provided that
{\roster
\itemitem{ $\otimes$ }  $|A|^n + |A| \le |{\Cal U}|$ or just $|R| \le
|{\Cal U}|$.
\endroster}
\endroster
\endproclaim
\bigskip

\demo{Proof}  Let $\{\langle a^j_\ell:\ell < n \rangle:j < |R|\}$ list the
$n$-tuples in $R$.  Choose $b_j \in {\Cal U} \backslash A$ for $j < |R|$ with 
no repetition.  For each $a \in A$ and $\ell < n$ let $Y^\ell_a = \{j:a^j_\ell
= a\}$, so clearly $a' \ne a'' \Rightarrow Y^\ell_{a'} \cap Y^\ell_{a''} =
\emptyset$ and let $\langle j_{a,\ell,k}:k < |Y^\ell_a| \rangle$ list
$Y^\ell_a$ with no repetition. \nl
Define $F^M_\ell$ by: $F^M_\ell(a) = b_{j_{a,\ell,0}},F^M_\ell(
b_{j_{a,\ell,k}}) = b_{j_{a,\ell,k+1}}$ except if $Y^\ell_a = \emptyset$
then $F^M_\ell(a_j) = a$.
\sn
Let

$$
\align
\varphi(\bar x,F_0,\dotsc,F_{n-1}) = (\exists z) \dsize \bigwedge_\ell 
[&F_\ell(x_\ell) \text{ well defined} \\
  &\and \neg(\exists y)(y \ne x_\ell \and F_\ell(y) = x_\ell) \and \theta
(x_\ell,z,F_\ell)]
\endalign
$$
\mn
where

$$
\align
\theta(x_\ell,z,F_\ell) =: \forall X&(x_\ell \in X \and (\forall y_1,y_2) \\
  &(y_1 \in X \and y_2 = F_\ell(y_1) \and y_1 \ne z \rightarrow y_2 \in
X) \rightarrow z \in X).
\endalign
$$
\mn
Those are monadic formulas.  Clearly,
\mr
\item "{$(*)_1$}"  $\varphi$ does not depend on ${\Cal U}$
\sn
\item "{$(*)_2$}"  $M \models \theta(a,b,F_\ell)$ \ub{iff} $b \in
\{F^{[i]}_\ell(a):i\}$ where $F^{[0]}_\ell(a) = a,F^{[\alpha]}_\ell(a) =
F_\ell(F^{[i]}_\ell(a))$ (if well defined)
\sn
\item "{$(*)_3$}"  $M \models ``F_\ell(a)$ well defined $\and \neg(\exists x)
(y \ne F(y) = a)$ \ub{iff} $a \in A$ [check].
\ermn
Hence if $\bar a = \langle a_\ell:\ell < n \rangle \in A$, by $(*)_1 +
(*)_2$ and definition of the $F^M_\ell$'s, $\varphi$ and $\theta$:

$$
\align
M \models \varphi[\bar a,F^M_0,\ldots] &\text{ \ub{iff} for some } z,
\dsize \bigwedge_{\ell < n} (a_\ell \in A \and z \in \{F^{(i)}_\ell(a):i\}) \\
  &\text{ \ub{iff} for some } j, \langle a_\ell:\ell < n \rangle = 
\langle a^j_\ell:\ell < n \rangle.
\endalign
$$
\sn
${{}}$  \hfill$\square_{\scite{1.12}}$
\enddemo
\bigskip

\demo{\stag{1.13} Conclusion}  If $\mu_1,\mu_2$ are functions with domain 
${\frak U},n < \omega$ and $(\forall {\Cal U} \in {\frak U})
[\mu_1({\Cal U})^n \le \mu_2({\Cal U})]$ \ub{then} 

$$
Q^{n\text{-ary}}_{\mu_1} \le_{\text{int}} Q^{1-1}_{\mu_2}
\text{ mod } Q^{\text{mon}}_{\le(\mu_1)^n}
$$
\mn
hence

$$
Q^{n\text{-ary}}_{\mu_1} \le_{\text{exp}} Q^{1-1}_{\mu_2}.
$$
\enddemo
\newpage

\head {\S2 Monadic analysis of $\exists_R$} \endhead  \resetall 
\bn
Our aim is to interpret $Q^{\text{mon}}_\lambda$ in $\exists_R$ for a 
maximal $\lambda$ and show that except on $\lambda$ elements $R$ is trivial.
So continuing later the analysis of $\exists_R$, we can instead analyze
$\{Q^{\text{mon}}_\lambda,\exists_{R_1}\}$ or analyze $\exists_{R_1}$ mod
$Q^{\text{mon}}_\lambda$ where $|\text{Dom}(R_1)| \le \lambda$ and 
$\exists_{R_1} \le_{\text{int}} Q_R$ and even $\exists_{R_1}
\le_{1\text{-int}} Q_R$ mod $Q^{\text{mon}}_\lambda$.  
This is made exact below.
\definition{\stag{2.1} Definition}  1) For any relation $R$ (on ${\Cal U}$)
let 

$$
\lambda_0 = \lambda_0(R) = \text{ Min}\{\frac{|{\Cal U}|}{2},\lambda'_0(R)\}
$$
\mn
where

$$
\align
\lambda'_0(R) = \text{ Min}\{|A|:&A \subseteq {\Cal U} \text{ and for every 
sequence } \bar b,\bar c \in {\Cal U} \\
  &(\text{of length } n(R)) \text{ we have } \bar b 
\approx_A \bar c \text{ implies } R[\bar b] \equiv R[\bar c]\}
\endalign
$$
\mn
where on $\approx_A$ see below  \nl
2)  $\bar b \approx_A \bar c$ means $\bar b = \langle b_i:i < n \rangle,
\bar c = \langle c_i:i < n \rangle$ and
\mr
\item "{$(a)$}"  $b_i \in A$ iff $c_i \in A$
\sn
\item "{$(b)$}"  $b_i \in A$ implies $b_i = c_i$
\sn
\item "{$(c)$}"  $b_i = b_j$ iff $c_i = c_j$.
\ermn
3) For a set $\Delta$ of formulas 
$\varphi(\bar x)$ (where $\varphi$ is a formula, 
$\bar x$ a finite sequence of variables including all variables occuring 
freely in $\varphi$) let

$$
\text{tp}_\Delta(\bar b,A,\bar M) = \{\varphi(\bar x,\bar a):\varphi
(\bar x,\bar y) \in \Delta,\bar a \subseteq A \text{ and }
M \models \varphi[\bar b,\bar a]\}.
$$
\mn
We omit $M$ when its identity is clear, and when $M = ({\Cal U},R)$ we
may write $R$ instead of $M$.  We may write ${\Cal U} \models \varphi
[\bar b,\bar a;R]$.  Replacing $\Delta$ by bs means $\Delta = \{\varphi
(\bar x):\varphi$ atomic or negation of atomic formula$\}$, here bs stands for
basic.  We may write
$\varphi$ instead $\{\varphi\}$ and $\Delta$ will be always finite. \nl
4) $S^m_\Delta(A,M) = \{\text{tp}_\Delta(\bar b,A,M):\bar b \subseteq M$
and $\ell g(\bar b) = m\}$.
\enddefinition
\bigskip

\remark{\stag{2.1A} Remark}  1) Note that 
$\lambda_0(R) \le \lambda'_0(R) \le |\text{Dom}(R)|$. \nl
2)  Note that if an equivalence relation $E$ on a subset of ${\Cal U}$
contains an equivalence class of cardinality $k \ge \frac{|{\Cal U}|}{2}$ 
or exactly $k \ge \frac{|{\Cal U}|}{2}$ singleton classes or
$k = |\text{Dom}(E)| \ge \frac{|{\Cal U}|}{2}$, then $\lambda'_0(R) = 
\lambda_0(R) = |{\Cal U}| - k < \frac{|{\Cal U}|}{2}$.  Otherwise, 
$\lambda'_0(R) > \frac{|{\Cal U}|}{2}$ 
and $\lambda_0(R) =  \frac{|{\Cal U}|}{2}$.
\endremark
\bn
The main result of this section is:
\proclaim{\stag{2.2} Theorem}  1)  $Q^{\text{mon}}_{\lambda_0(R)}
\le_{\text{int}} \exists_R$; we mean, of course, uniformly. \nl
2)  There is a relation $R_1$ on ${\Cal U}$ with $n(R_1) = n(R)$ and,
$|\text{Dom}(R_1)| \le \lambda'_0(R) + n$ such that $\exists_R 
\equiv_{\text{int}} \{\exists_{R_1},Q^{\text{mon}}_{\lambda'_0(R)}\}$.  In
fact, $\exists_R \equiv_{1-\text{int}} \exists R_1 \text{ mod }
Q^{\text{mon}}_{\lambda_0(R)}$. 
\endproclaim
\bn
The proof is broken into some claims.
\proclaim{\stag{2.3} Claim}  Let $R$ be an $n$-place relation on ${\Cal U}$
such that $n > 1$.
We can find a set $A$, sequences $\bar a_i$ and
elements $b_i,c_i$ for $i < i^*$, where $i^* \ge
\frac{\lambda'_0(R)}{n(R)(n(R)-1)}$ or
$i^* \ge \frac{|{\Cal U}| - n(R)}{n(R)(n(R)-1)}$ such that:
\mr
\item "{$(a)$}"   $\bar a_i$ is with no repetition and is $\subseteq A$
\sn
\item "{$(b)$}"  $b_i,c_i \notin A$,
\sn
\item "{$(c)$}"  $\langle (b_i,c_i):i < i^* \rangle$ is 
with no repetition, i.e. $i \ne j \Rightarrow b_i \ne b_j
\and c_i \ne c_j \and b_i \ne c_j$
\sn
\item "{$(d)$}"  tp$_{\text{bs}}(\bar a_i \char 94 \langle b_i \rangle,
\emptyset,R) \ne \text{ tp}_{\text{bs}} (\bar a_i \char 94 \langle c_i 
\rangle,\emptyset,R)$, that is for some atomic formula 
$\varphi$ we have $({\Cal U},R) \models \varphi (\bar a_i,b_i) \equiv 
\neg \varphi(\bar a_i,c_i)$.
\endroster
\endproclaim
\bigskip

\demo{Proof}  We try to choose by induction on $i,\langle A^i_\ell:\ell <
n(R) \rangle,\langle \bar a_i,b_i,c_i \rangle$ and $\ell(i) < n(R)$ such
that:
\mr
\widestnumber\item{$(viii)$}
\item "{$(i)$}"  $\ell < k < n(R) \Rightarrow A^i_\ell \cap
A^i_k = \emptyset$
\sn
\item "{$(ii)$}"  $A^i_\ell \subseteq A^i_{\ell +1}$ and 
$A^0_\ell = \emptyset$
\sn
\item "{$(iii)$}"  $\bar a_i \char 94 \langle b_i \rangle \char 94
\langle c_i \rangle$ is with no repetition and has length $\le n(R)+1$
\sn
\item "{$(iv)$}"  tp$_{\text{bs}}(\bar a_i \char 94 \langle b_i \rangle,
\emptyset,R) \ne \text{ tp}_{\text{bs}}(\bar a_i \char 94 \langle c_i 
\rangle,\emptyset,R)$, that is for some atomic formula $\varphi(\bar x,y)$
(so gotten from $R(x_0,\dotsc,x_{n(R)-1})$ by substitution) we have \nl
$\varphi(\bar a_i,b_i) \equiv \neg \varphi(\bar a_i,c_i)$
\sn
\item "{$(v)$}"  $b_i,c_i \notin \cup\{A^i_\ell:\ell < n(R)\}$
\sn
\item "{$(vi)$}"  $\ell(i) = \text{ Min}\{\ell:\bar a_i \cap A^i_\ell 
= \emptyset\}$
\sn
\item "{$(vii)$}"  $A^{i+1}_{\ell(i)} = A^i_{\ell(i)} \cup \{b_i,c_i\}$
\sn
\item "{$(viii)$}"  $A^{i+1}_\ell$ is $A^i_\ell \cup \{\text{Rang}(\bar a_i)
\backslash \dbcu_m A^i_m\}$ \nl
if $(\ell = 0 \and \ell(i) > 0) \vee (\ell =1 \and \ell(i) = 0)$
\sn
\item "{$(ix)$}"  $A^{i +1}_\ell = A^i_\ell$ in the other cases.
\ermn
So for some $i = i(*)$ we cannot continue; we claim that 
$A =: \dbcu_\ell A^i_\ell$ has cardinality $\ge \lambda'_0(R)$ or
$\ge |{\Cal U}| - n(R)$. \nl
Why?  Otherwise by the definition of $\lambda'_0(R)$ there are sequence
$\bar b,\bar c$ from ${\Cal U}$ of length $n(R)$ such that $\bar b \approx_A
\bar c$ but $\bar b \in R \equiv \bar c \notin R$.  Hence we can 
find sequences $\bar b',\bar c'$ from ${\Cal U}$ of the same length 
$\le n(R)$, each with no
repetitions such that $\bar b' \approx_A \bar c'$ but for some $\varphi =
R(\bar x) = R(x_{i_0},x_{i_1},\dotsc,x_{i_{n_R}-1}),\ell g(\bar x) = \ell g
(\bar b') = \ell g(\bar c')$ we have $\varphi(\bar b') \and \neg \varphi
(\bar c')$.  Now we can find $k$ and $\bar d_0,\dotsc,\bar d_k$ such that:
$\bar d_0 = \bar b',\bar d_k = \bar c'$, and $\bar d_\ell$ is with no
repetitions, $\ell < k \Rightarrow \bar d_\ell \approx_A \bar d_{\ell +1}$
and $\ell < k \Rightarrow (\exists^!\, i)d_{\ell,i} \ne d_{\ell +1,i}$; here
we use the assumption toward contradiction $|A| \le |{\Cal U}|-n$.  
\sn
So for some $\ell < k$ we have $\varphi(\bar d_\ell) \and 
\neg \varphi(\bar d_{\ell +1})$.
Now let $r$ be such that $d_{\ell,r} \ne d_{\ell +1,r}$, so \wilog \,
$r = \ell g(\bar x)-1$, let $\bar a_{i(*)} = \bar d_{\ell,r} \restriction
(\ell g(\bar b_\ell)-1),b_{i(*)} = d_{\ell,r},c_{i(*)} = d_{\ell +1,r}$.
Clearly they are as required in clause (iii) + (iv), now $\ell(i(*))$ is 
well defined by clause (vi) as $|\text{Rang}(\bar a_i)| \le n(R)-1$, so 
$|\text{Rang}(\bar a_i) \cap (\dbcu_\ell A^i_\ell)| < n(R)$.
Now we can define $A^{i(*)+1}_\ell$ for $\ell < n(R)$ by clauses (vii),
(viii) and (ix).  Trivially, clauses (i) and (ii) hold, and we get a
contradiction to the choice of $i(*)$.  So really $|A| \ge \lambda'_0(R)$
or $|A| \ge |{\Cal U}| - n(R)$.
\mn
Now note that $|\dbcu_\ell A^{i(*)}_\ell| \le (n(R) +1) \times
i(*)$, by clauses (vii), (viii), (ix) so $i(*) \ge |\dbcu_\ell A^{i(*)}
_\ell|/(n(R)+1)$ and for some $\ell$ we have
$|\{i < i(*):\ell(i) = \ell\}| \ge i(*)/n(R)$.
So if $\lambda'_0(R) < |{\Cal U}| - n(R)$ we get $|\{i < i(*):\ell(i) =
\ell\}| \ge i(*)/n(R) \ge \lambda'_0(R)/((n(R)+1)(n(R))$.
If $\lambda'_0(R) \ge |{\Cal U}| - n(R)$ we get $|\{i < i(*):\ell(i) =
\ell\}| \ge (|{\Cal U}|-n(R))/(n(R)(n(R)+1))$.  So renaming we are done.
\hfill$\square_{\scite{2.3}}$
\enddemo
\bigskip

\proclaim{\stag{2.4} Claim}  There is a formula $\varphi^* = \varphi^*
(x,\bar y;R)$, in first order logic, of course, such that:
\mr
\item "{$(*)$}"  if ($R$ is an $n(R)$-place relation on ${\Cal U}$ and)
$\lambda'_0(R) < \frac 23|{\Cal U}|$, \ub{then} $\varphi^*$ exemplifies
$Q_{\lambda_0(R)} \le_{\text{int}} Q_R$ even $Q_{\lambda_0(R)}
\le_{1\text{-int}} Q_R$ specifically, for some \nl
$\bar d$ we have $\{a:({\Cal U},R) \models \varphi^*(a,\bar d,R)\}$ has 
$\lambda_0(R)$ members.
\endroster
\endproclaim
\bigskip

\demo{Proof}  Without loss of generality $\frac 13|{\Cal U}| > n(R)^2 +n(R)$.

Let $A \subseteq {\Cal U}$ be a set of power $\lambda'_0(R)$ such that $\bar b
\approx_A \bar c$ implies $R[\bar b] \equiv R[\bar c]$.  As $|{\Cal U}| -
\lambda'_0(R)$ is large enough, we can find pairwise distinct 
$d_i \in {\Cal U} \backslash A$ for $i < n(R)^2$.  Define $\bar d = 
\langle d_i:i < n(R)^2 \rangle$ and $\varphi^*(x,\bar d,R) = 
\bigvee\{(\exists y_0,\dotsc,y_{k-1})$ 
[the elements $y_0,\dotsc,y_{k-1},x$ are pairwise distinct and for any $m$ if
the elements $y_0,\dotsc,y_{k-1},d_m,x$ are pairwise distinct and 
$\varphi(x,y_0,\dotsc,y_{k-1}) \equiv \neg \varphi(d_m,y_0,\dotsc,y_{k-1})]:
\varphi = \varphi
(z_0,\dotsc,z_k,R)$ is an atomic formula in $L(R)$ (so $k+1 \le n(R)$) and
$m < n(R)^2$, so $m,k$ are natural numbers$\}$.  By the choice of $A$ we have
$x \notin A \Rightarrow \neg \varphi^*(x,\bar d,R)$, hence $B =: \{x \in 
{\Cal U}:{\Cal U} \models \varphi^*[x,\bar d,R]\}$ is a subset of $A$.  
Clearly
$Q^{\text{mon}}_{|B|} \le_{\text{int}} \exists_R$ (uniformly); hence it
suffices to prove $|B| = \lambda'_0(R)$ which follows if we show
\mr
\item "{$(*)$}"  if $\bar b \cong_B \bar c$ then $R[\bar b] \equiv R[\bar c]$.
\ermn
For this it suffices to prove
\mr
\item "{$(**)$}"  if $\varphi(\bar x,R) \in L(R)$ is atomic, $\bar b,\bar c$
are sequences of length $\ell g(\bar x) \le n(R)$ without repetition then
$\bar b \cong_B \bar c$ implies $\varphi(\bar b,R) \equiv \varphi(\bar c,R)$.
\ermn
To prove $(**)$, by reordering the sequences we let $\bar b \char 94 
\bar c_0,\bar b \char 94 \bar c_1$ be sequences from
${\Cal U}$, without repetition, $\bar b \subseteq B,\bar c_0,\bar c_1$
disjoint to $B$; by the transitivity of $\equiv$, without loss of
generality $\bar c_1$ is disjoint to $\bar d$.  Now for some 
$i,\langle d_i,d_{i+1},\dotsc,d_{i+k-1}
\rangle$ (where $k = \ell(\bar c_0))$ is disjoint to $\bar c_0$ (and
obviously to $\bar c_1$).
\sn
Now we shall prove that for every atomic $\varphi(\bar x,\bar y,R),\ell g
(\bar x) = k,\ell g(\bar y) = \ell g(\bar b)$ we have 
$\models \varphi(\bar c_\ell,\bar b,R)
\equiv \varphi (\langle d_i,\dotsc,d_{i+k-1}\rangle,\bar b,R)$ thus finishing.
For this we define $\bar c_{\ell,m}(m \le k)$ such that each $\bar c_{\ell,m}$
is with no repetitions, disjoint to $B \cup \bar b$ and $\bar c_{\ell,0} = 
\bar c_\ell,\bar c_{\ell,k} = \langle d_i,\dotsc,d_{i+k-1} \rangle$,
$\bar c_{\ell,m+1},\bar c_{\ell,m}$ are distinct in one place only.  
By the definition of $B$ (and $\varphi$) for every atomic
$\varphi(\bar x,\bar y,R)$ we have $\models \varphi(\bar c_{\ell,m},\bar b,R) 
\equiv \varphi (\bar c_{\ell,m+1},\bar b,\bar R)$ so we finish easily.
(Being more careful, e.g. $\frac 19 |{\Cal U}| \ge n(R)$ suffices). \nl
${{}}$  \hfill$\square_{\scite{2.4}}$
\enddemo
\bigskip

\remark{Remark}  Note that definition of $\lambda_0$ applies to any
relation, in particular, the relation being defined by a formula so we 
may freely speak at $\lambda_0(\psi)$ or $\lambda_0(\psi(\bar x))$.
\endremark
\bigskip

\proclaim{\stag{2.5} Claim}  
$Q^{\text{mon}}_{\le \lambda_0(R)} \le_{\text{int}} Q_R$.
\endproclaim
\bigskip

\demo{Proof}  If we can replace in $R$ some variables by constants or
other variables having at least one equality getting a relation $R'$ such
that $\lambda_0(R') \ge \text{Min}\{\lambda_0(R),\frac{1}{7n(R)}|{\Cal U}|\}$ 
we do it: or in other words we are inducting on $n(R) \ge 1$.  
\bn
\ub{Case 1}:  $n(R) = 1$.

So $R$ is unary; now note that each of the sets $A = R,A' = {\Cal U}
\backslash R$ can serve in the definition of $\lambda'_0(R)$, hence

$$
\lambda'_0(R) \le \text{ Min}\{|A|,|A'|\} = \text{ Min}\{R,|{\Cal U}
\backslash R|\} \le \frac{|{\Cal U}|}{2}
$$
\mn
so we are clearly done.
\bn
\ub{Case 2}:  $n(R) > 1$.
\mn
If $\lambda'_0(R) \le \frac 23 |{\Cal U}|$ we can interpret $Q^{\text{mon}}
_{\le \lambda_0(R)}$ by \scite{2.4}, as it suffices to show that at least one
of several $\varphi$'s interpret.  So assume $\lambda'_0(R) > \frac 23 
|{\Cal U}|$.  
Hence $\lambda_0(R) = \frac{|{\Cal U}|}{2}$ and we shall prove that we can
interpret $Q^{\text{mon}}_{\le \lambda_0(R)}$, for this it is enough if we
can show that we can interpret 
$Q^{\text{mon}}_{\le[\frac{|{\Cal U}|}{2^{2^{n(R)}}}]}$.  For this is enough
to find first order $\theta(\bar x_1,\bar y_1,R),\dotsc,\theta_k(\bar x_k,
\bar y_k,R)$ with the $k$ and $\theta_\ell$ depending only on $n(R)$ and
not on $|{\Cal U}|$ such that $\ell g(\bar x_\ell) < n(R)$ and for some
$\ell \in \{1,\dotsc,k\}$ and $\bar b \in {}^{\ell g(\bar y_\ell)}{\Cal U}$,
we have

$$
\frac{|{\Cal U}|}{2^{2^{n(R)}}} \le \lambda'_0(\theta_\ell(-,\bar b,R)).
$$
\mn
For any $\ell < k < n(R)$ we can consider the formula
$R_{\ell,k}(x_0,\dotsc,x_{n(R)-1}) = R(x_0,\dotsc,x_{n(R)-1}) \and
x_\ell = x_k$ and $R^* = R(x_0,\dotsc,x_{n-1}) \and
\dsize \bigwedge_{\ell < k} x_\ell \ne x_k$.  Easily $\lambda'_0(R) \le
\lambda'_0(R^*) + \dsize \sum_{\ell < k < n(R)} \lambda'_0(R_{\ell,k})$.
Now if for some $\ell < k,\lambda_0(R_{\ell,k}) \ge
\frac{|{\Cal U}|}{3n(R)(n(R-1))}$ we are done by the induction hypothesis.
So we can assume $\lambda'_0(R^*) \ge \frac{|{\Cal U}|}{3n(R)(n(R)-1)}$
hence by the above \wilog \, $\lambda'_0(R^*) \ge \frac 23 |{\Cal U}|$, so
we can assume 
\mr
\item "{$(*)_0$}"  $R = R^*$.
\ermn
So atomic formulas not equivalent to a fix truth value except 
equality are just $R(\ldots,x_{\sigma(\ell)},\ldots)$ for 
$\sigma \in \text{ Per}(n(R))$. \nl
Let $A,\bar a_i,b_i,c_i$ for $i < j^* = (|{\Cal U}| - n(R))/(n(R)(n(R)-1))$ 
be as guaranteed by \scite{2.3}.
For some atomic $\varphi = \varphi(\bar x,y) = \varphi(\bar x,y,R)$ we have 
$|\{i:\varphi(\bar a_i,b_i) \wedge \neg \varphi(\bar a_i,c_i)\}| \ge 
j^*/n(R)$, by $(*)_0$.  Without loss of generality this occurs for 
$i < j^*/n(R)$.  For $i \le j^*$ let $F_i$ be the permutation of 
${\Cal U}$, interchanging $b_j,c_j$ for $j<i$ and being the identity 
otherwise.  Let $R_i = F''_i(R),\psi_i = \psi_i(\bar x,y) =: [\psi(\bar x,y,R,
R_i) = \varphi(\bar x,y,R) \and \neg \varphi(\bar x,y,R_i)]$ so 
\mr
\item "{$(*)$}"  $\psi_j(\bar a_i,b_i) \and \neg \psi_j(\bar a_i,c_i)$ 
if $i<j$.
\ermn
So by the definition of $\lambda'_0(-)$ we have $\lambda'_0(\psi_j) \ge j$ 
for $j \le j^*/n(R)$, where we consider $\psi_j$ as a
$(\ell g(\bar x)+1)$-place relation. \nl
[Why?  If $A \subseteq |{\Cal U}|,|A| < j$ exemplifies the failure of this
assertion (by the definition of $\lambda'_0(R)$) then $w =: \{i < j:A \cap
\{b_i,c_i\} \ne \emptyset\}$ has $\le |A|$ members, so choose $i \in j
\backslash w$, now ${\Cal U} \models \psi_j(\bar a_i,b_i,c_i) \and \neg
\psi_j(\bar a_i,c_i,b_i)$, (holds by $(*)$) contradict the choice of $A$.]
\nl
So if $\dsize \bigvee_j [|{\Cal U}|/3 \le \lambda'_0(\psi_j) < 
\frac{2}{3}|{\Cal U}|]$ we are done; hence assume not.
\mn
If for every $j$ we have $\lambda'_0(\psi_j) < 
\frac{|{\Cal U}|}{3}$, then we get
$[\frac{|{\Cal U}|-n(R)}{n(R)^2(n(R)-1)}] = [j^*/n(R)] \le \lambda'_0
(\psi_{[j^*/n(R)]}) \le \frac{|{\Cal U}|}{3}$, so we easily finish by
\scite{2.4}. \nl
Also $\lambda'_0(\psi_0) = \lambda'_0(\emptyset)=0$ as $R_0=R$ so
${\Cal U} \models \varphi(\bar x,y,R) \equiv \varphi(\bar x,y,R_0)$ hence
${\Cal U} \models \neg \psi_0(\bar x,y,R,R_0)$.  Without loss of generality
$\varphi$ is $R$.    So the bad case is
that for some $j$ we have $\lambda'_0(\psi_j) < \frac 13|{\Cal U}|$ and 
$\lambda'_0(\psi_{j+1}) \ge \frac 23 |{\Cal U}|$.
Let $B^* \subseteq {\Cal U}$ exemplify $\lambda'_0(\psi_j) < \frac 13
|{\Cal U}|$.  Let for $\ell < n(R),
\theta_\ell(x_0,\dotsc,x_{n(R)-2},y,R) = R(x_0,\dotsc,x_{\ell-1},y,
x_\ell,\dotsc,x_{n(R)-2})$ and on $\theta_\ell(\langle x_m:m < n(R)-1 
\rangle,c_j;R_j)$ we apply our induction
hypothesis as its arity is $\ell g(\bar x)$ which is at most $n(R)-1$ (see the
beginning of the proof) hence 
$\lambda'_0(\theta_\ell(\langle x_m:m < n(R)-1 \rangle,c_j,R)) \le
\frac{1}{7n(R)}|{\Cal U}|$ and let $B_\ell \subseteq {\Cal U}$ exemplify it.
Similarly let $B'_\ell \subseteq {\Cal U}$ exemplify $\lambda'_0(\theta_\ell
(\langle x_m:m < n(R)-1 \rangle,b_j,R)) \le \frac{1}{7n(R)}$.
Let $B = B^* \cup \dbcu_{\ell < n(R)} B_\ell \cup \dbcu_{\ell < n(R)}
B'_\ell \cup \{b_j,c_j\}$.  Now $B$ is a subset of
${\Cal U}$ with $< (\frac{2}{7} + \frac 13)|{\Cal U}| + 2 < 
\frac 23|{\Cal U}|$ elements.  By the definition of $\psi_j,\psi_{j+1}$ 
and $(*)_0$ such $B$ exemplifies $\lambda'_0(\psi_{j+1}) < 
\frac 23|{\Cal U}|$, contradiction.
\nl
${{}}$   \hfill$\square_{\scite{2.5}}$
\enddemo
\bn
We have implicitly used:
\proclaim{\stag{2.6} Claim}  If $R$ is a Boolean combination of 
$R_0,\dotsc,R_{n-1}$ \ub{then} $\lambda'_0(R) \le
\dsize \sum_{\ell < n} \lambda'_0(R_\ell)$ hence 
$\lambda_0(R) \le \dsize \sum_{\ell < n} \lambda_0(R_\ell)$.
\endproclaim
\bigskip

\demo{Proof}  If $A_\ell$ witnesses the value $\lambda'_0(R_\ell)$ then
$A = \dbcu_{\ell < k} A_\ell$ witnesses $\lambda'_0(R) \le |A| \le
\dsize \sum_{i < k} |A_\ell|$.  \hfill$\square_{\scite{2.6}}$
\enddemo
\bn
Now we turn to \scite{2.2}
\demo{Proof of \scite{2.2}(1)}

Immediate by \scite{2.4}, \scite{2.5}.
\enddemo
\bn
\centerline{$* \qquad * \qquad *$}
\bigskip

\demo{Proof of \scite{2.2}(2)}  
Let $d_i$ (for $i < n(R))$ be distinct elements of ${\Cal U} \backslash A$
where $A$ exemplifies $\lambda'_0(R)$ as if $\lambda'_0(R) +n \ge
|{\Cal U}|$ then we can choose $R_1 = R$.  
Of course, we can concentrate on the case $n(R) > 1$.  Let
$R_1 = R \restriction (A \cup \{d_i:i < n(R)\}$.
So $\langle a_1,\dotsc,a_n \rangle \in R$ \ub{iff} for some $\langle
a'_1,\dotsc,a'_n \rangle \in R_1$ we have $\langle a_1,\dotsc,a_n \rangle
\approx_A
\langle a'_1,\dotsc,a'_n \rangle$ and $\dsize \bigwedge_\ell[a'_\ell \notin
A \rightarrow \dsize \bigvee_m[a'_\ell = d_m]$, so we can define $R_1$ from
$R$ and $R$ from $R_1$ by a
quantifier free formula using the unary relation $A$ and individual constants
$d_0,d_1,\dotsc,d_{n(R)-1}$.  Hence $\exists_R \le_{1\text{-int}}
\exists_{R_1}$ mod $Q^{\text{mon}}_{\lambda'_0(R)}$ but $\lambda'_0(R) \le
2 \lambda_0(R)$ so $\exists_R \le_{1\text{-int}} \exists_{R_1}$ mod
$Q^{\text{mon}}_{\lambda_0(R)}$.

Also easily $\{\exists_{R_1},Q^{\text{mon}}_{\lambda_0(R)}\} 
\le_{1-\text{int}} \exists_R$.  \hfill$\square_{\scite{2.2}}$
\enddemo
\bn
We can get the parallel result for $Q_K$.
\definition{\stag{2.7} Definition}  Let $\lambda_0(K) = \text{ Min}\{\lambda:
R \in K \Rightarrow \lambda_0(R) < \lambda\}$ note that the minimum is taken
for each ${\Cal U} \in {\frak U}$ separately.
\enddefinition
\bigskip

\proclaim{\stag{2.8} Theorem}  1) $Q^{\text{mon}}_{\le \lambda_0(K)}
\le_{\text{int}} \exists_K$. \nl
2)  There is $K_1,n(K_1) = n(K)$ such that
\mr
\item "{$(a)$}"  $\exists_K \equiv_{\text{int}} \{\exists_{K_1},
Q^{\text{mon}}_{< \lambda_0(K)}\}$
\sn
\item "{$(b)$}"  $R \in K \Rightarrow |\text{Dom}(K)| < \lambda_0(K)$
\sn
\item "{$(c)$}"  $\exists_K \equiv_{1\text{-int}} \exists_{K_1}$ mod
$Q^{\text{mon}}_{< \lambda_0(K)}$.
\endroster
\endproclaim
\bigskip

\demo{Proof}  Immediate by the uniformity of our results.
\enddemo
\bn
\ub{\stag{2.9} Discussion}:  The interpretation here uses first order formulas
of low complexity but use several copies of $R$.  We may wonder if we can
just use one copy of $R$ by complicating the formula.  Now if $R$ is a
connected graph every node having a valency $\le m << |{\Cal U}|$, we see that
not.  But we can prove that the general situation in the problematic case 
is not far from this (similar to a model of a strongly minimal theory, a local
version).  Also in general 2 copies of $R$ suffice.
\newpage

\head {\S3 The one-to-one function analysis} \endhead  \resetall 
\bn
The aim of this section is similar to the previous one, going one step
further, i.e. we want to analyze $\exists_R$, interpreting in it 
$Q^{1-1}_\lambda$ for a maximal $\lambda$, hoping that ``the remainder" has
domain $\le \lambda$.
\definition{\stag{3.1} Definition}  Let $\lambda_1(R)$ be
Max$\{|\{\text{tp}_{\text{bs}}(a,A,R):a \in {\Cal U} \backslash A\}|:
A \subseteq {\Cal U}\}$.  (On tp$_{\text{bs}}$ see \scite{2.1}(3)).
\enddefinition
\bigskip

\demo{\stag{3.2} Fact}  $\lambda_1(R) \le \lambda'_0(R)+1$ and if equality
holds then $\lambda_1(R) \le 2^{2^{n(R)^2}}$.
\enddemo
\bigskip

\demo{Proof}  Straight, assume $A_0$ exemplifies $\lambda'_0(R)$ and let
$A \subseteq {\Cal U}$.  Then $a,b \in ({\Cal U} \backslash A) \backslash 
A_0 \Rightarrow \text{tp}_{\text{bs}}(a,A,R) = \text{ tp}_{\text{bs}}(b,A,R)$ 
by the choice of $A_0$ hence
$|\{\text{tp}_{\text{bs}}(a,A,R):a \in {\Cal U} \backslash A\}| \le |A_0
\backslash A| + 1 \le |A_0| + 1 = \lambda'_0(R)+1$.  Next assume that equality
holds, so necessarily $|A_0 \backslash A| = |A_0|$ hence $A \cap A_0 = 
\emptyset$; now choose $A' \subseteq A$ with Min$\{n(R)-1,|A|\}$ elements.  
By the
choice of $A_0$, if $b,c \in {\Cal U} \backslash A$ then

$$
\text{tp}_{\text{bs}}(b,A,R) = \text{ tp}_{\text{bs}}(c,A,R) \Leftrightarrow
\text{ tp}_{\text{bs}}(b,A',R) = \text{ tp}_{\text{bs}}(c,A',R).
$$
\mn
[Why?  $\Rightarrow$ holds as $A' \subseteq A$; next we shall prove 
$\Rightarrow$.  This suffices so assume tp$_{\text{bs}}(b,A',R) =
\text{ tp}_{\text{bs}}(c,A',R)$.  So let $\varphi(x,\bar y,R)$ be an atomic
formula (i.e. a substitution in $R(x_0,\dotsc,x_{n(R)-j})$, so
$\ell g(\bar y) + 1 \le n(R))$ and let $\bar a_1$ be a sequence of length
$\ell g(\bar y)$ from $A$, we shall show that $\varphi(b,\bar a_1,R) \equiv
\varphi(c,\bar a_1,R)$, this suffices.  If $|A| < n(R)$, then $A'=A$ and we
are done, so assume $|A| \ge n(R)$.

We can find a sequence $\bar a_2$ from $A'$ which realizes the same equality
type as $\bar a_1$ (because $\ell g(\bar a_1) = \ell g(\bar y) \le n(R)-1 =
|A'|$).  Now by our assumption $\varphi(b,\bar a_2,R) \equiv \varphi(c,
\bar a_2,R)$ (that is as tp$_{\text{bs}}(b,A',R) = \text{ tp}_{\text{bs}}(c,
A',R))$, so to get our desired $\varphi(b,\bar a_1,R) \equiv \varphi(c,
\bar a_1,R)$ it suffices to prove $\varphi(b,\bar a_1,R) \equiv \varphi(b,
\bar a_2,R)$ and $\varphi(c,\bar a_1,R) \equiv \varphi(c,\bar a_2,R)$.  But
on both $b$ and $c$ we just assume they are in ${\Cal U} \backslash A$, so
by symmetry it is enough to show $\varphi(b,\bar a_1,R) \equiv \varphi(b,
\bar a_2,R)$.  Now as $\bar a_1,\bar a_2$ are included in $A$ and
have the same equality type (over the $\emptyset$), by the choice of $A_0$
and as $A_0 \cap A = \emptyset$ necessarily $\bar a_1,\bar a_2$ realizes the
same equality type over ${\Cal U} \backslash A$, so as $b \in {\Cal U}
\backslash A$ we have $\varphi(b,\bar a_1,R) \equiv \varphi(b,\bar a_2,R)$.]
\nl
Hence
$\lambda_1(R) \le |\{\text{tp}_{\text{bs}}(b,A',R):b \in {\Cal U}\}| \le
2^{|\Phi|}$ where $\Phi$ is the set of atomic formulas $\varphi(x,\bar a)$
such that $\bar a \subseteq A',|\Phi| \le n(R) \times (n(R)-1)^{n(R)-1} \le
2^{n(R)^2}$.  \hfill$\square_{\scite{3.2}}$
\enddemo
\bigskip

\proclaim{\stag{3.3} Claim}  $Q^{1-1}_{\lambda_1(R)} \le_{\text{int}}
\exists_R$; of course uniformly.
\endproclaim
\bigskip

\demo{Proof}  Suppose $h$ is a one-to-one, one place partial function from
${\Cal U}$ to ${\Cal U}$ with $\lambda = |\text{Dom}(h)| \le \lambda_1(R)$ and
$\lambda \le \frac{1}{n(R)+1}|{\Cal U}|$ (we use freely \scite{1.3}).  Let $A
\subseteq {\Cal U}$ be such that $\{\text{tp}_{\text{bs}}(a,A,R):a \in 
{\Cal U} \backslash A\}$ has cardinality $\lambda_1(R)$.  So we can
find $a_i \in {\Cal U} \backslash A$ (for $i < \lambda$) such that
$\text{tp}_{\text{bs}}(a_i,A,R)$ are pairwise distinct.  Retaining the last
sentence (by not necessarily the original demand on $A$) without loss 
of generality $|A| \le |{\Cal U}| - \lambda - \lambda$. \nl
[Why?  Just for each $i=1,\dotsc,\lambda-1$ choose $\bar d_i \subseteq A$ of
length $< n(R)$ such that $\langle \text{tp}_{\text{bs}}(a_j,
\dbcu^i_{\ell=1} \bar d_\ell):j \le i \rangle$ is with no repetitions so
\wilog \, $|A| \le (n(R)-1) \times \lambda$ and compute.]
Let $h = \{\langle b_i,c_i \rangle:i < \lambda\}$, without loss of
generality $b_i,c_i \notin A$ (just permute $R$, i.e. using an 
isomorphic $R'$) and we can find $F_1,F_2$
permutation of ${\Cal U}$ which are the identity on $A$ such that
$F_1(a_i) = b_i,F_2(a_i)=c_i$.
Let $R_1 = F_1(R)$ and $R_2 = F_2(R)$ and define the monadic relations
$P_0 = A,P_1 = \{b_i:i < \lambda\},P_2 = \{c_i:i < \lambda\}$ (all of
cardinality $\le \lambda_0(R))$.  
Let $\varphi(x,y,P_0,P_1,P_2,R_1,R_2)$ ``say" that for
every atomic $\psi(x,\bar z,R) \in L(R)$ and $\bar t \in P_0$ we have:
$\varphi(x,\bar t,R_1) \equiv \varphi(y,\bar t,R_2)$ and $P_1(x),P_2(y)$.  
Clearly $\varphi$ defines $h$.  \hfill$\square_{\scite{3.3}}$
\enddemo
\bigskip

\proclaim{\stag{3.4} Lemma}  Assume $\lambda_1(R) \times n(R)^2 +
n(R) < |{\Cal U}|$.  For any set $A \subseteq {\Cal U}$, let $E_A$ be 
the following equivalence relation on ${\Cal U}$: tp$_{\text{bs}}(a,A,R) = 
\text{ tp}_{\text{bs}}(b,A,R)$.
For any $A \subseteq {\Cal U}$ and $\bar C = \langle C_\ell:\ell < k \rangle$
such that $C_\ell \subseteq {\Cal U}$ let $E_{A,{\bar C}}$ be the following 
equivalent relation on ${\Cal U}: a E_{A,{\bar C}} b$ iff $a E_A b \and
\dsize \bigwedge_\ell a \in C_\ell \equiv b \in C_\ell$.  There are a set 
$A \subseteq {\Cal U}$ and sequence $\bar C = \langle C_\ell:\ell < n(R)-2
\rangle$ with $C_\ell \subseteq {\Cal U}$ such that
\mr
\item "{$(A)$}"  $|A| \le n(R) \times n(R) \times \lambda_1(R)$
\sn
\item "{$(B)$}"   if $\bar b \cong_\emptyset \bar c$ and 
$b_i E_{A,{\bar C}} c_i$ for all 
$i < \ell g(\bar b)$ \ub{then} $R(\bar b) \equiv R(\bar c)$
\sn
\item "{$(C)$}"  $E_A$ has at most $|A| + \lambda_1(R)$ classes
\sn
\item "{$(D)$}"  each $C_\ell$ has at most $\lambda_1(R)$ elements.
\endroster
\endproclaim

\demo{Proof}  We try by induction on $i$ to choose $\langle A^i_\ell:\ell <
n(R) \rangle$ such that
\mr
\widestnumber\item{$(iii)$}
\item "{$(i)$}"  $A^i_\ell \subseteq {\Cal U}$
\sn
\item "{$(ii)$}"  $\ell < k < n(R) \Rightarrow A^i_\ell \cap A^i_k =
\emptyset$
\sn
\item "{$(iii)$}"  $|A^i_\ell| \le i$ 
\sn
\item "{$(iv)$}"  $\dsize \sum_{k < n(R)} |\{\text{tp}_{\text{bs}}(b,
\dbcu_{\ell \ne k} A^i_\ell):b \in A^i_k\}|$ is at least $i$
\sn
\item "{$(v)$}"  $j < i \Rightarrow A^j_\ell \subseteq A^i_\ell$.
\ermn
Now for $i=0$ let $A^i_\ell = \emptyset$. \nl
We necessarily are stuck for some $i = i(*) \le
\lambda_1(R) \times n(R)$; i.e. $A^j_\ell$ are defined for $j \le i(*)$
but we cannot choose $\langle A^{i(*)+1}_\ell:\ell < n(R) \rangle$,
otherwise by clause (iv) for some $k$ the set $\{\text{tp}_{\text{bs}}(b,
\dbcu_{\ell \ne k} A^i_\ell):b \in A^i_k\}$ has at least $i/n(R)$ elements 
which is (by the assumption toward contradiction) $> \lambda_1(R)$,
but now $A = \dbcu_{\ell \ne k} A^i_\ell$ contradicts the definition of
$\lambda_1(R)$ as $A^i_k \cap A = \emptyset$ by clause (ii).  Let
$A = \dbcu_{\ell < n(R)}A^{i(*)}_\ell$.
For $\ell < n(R)-2$, choose $C_\ell$ as a set of representatives for $\{a/E_A:
a/E_A$ has $\le n(R)$ but at least $2 + \ell$ elements$\}$, such that
$C_\ell$ is disjoint to $\dbcu_{m < \ell} C_m$ and we shall show 
that $A,\bar C$ is as
required. Now clause $(A)$ holds by clause (iii) and the choice of $A$ (and
the bound above on $i(*)$).  Toward proving clause (B) assume 
$\bar b \cong_\emptyset \bar c$ and $b_\ell E_{A,{\bar C}} c_\ell$ for 
$\ell < \ell g(\bar b)$.  Without loss of generality
$\bar b$ has no repetitions.  
Note if $b_\ell \ne c_\ell$ then $b_\ell,c_\ell \notin A$ (as $\bar b
\approx_A c$) and $b_\ell E_A c_\ell$ (see definition of $E_{A,{\bar C}}$), 
but by the choice of the $C_\ell$'s, $b_\ell
/E_A = c_\ell/E_A$ has $> n(R)$ elements, so there is $d \in b_\ell/E_A
\backslash \{b_k:k < \ell g(\bar b)\}$.  Hence by transitivity of all 
the relevant
conditions \wilog \, for some $k(*) < \ell g(\bar b)$ we have
$b_{k(*)} \ne c_{k(*)} \and \dsize \bigwedge_{m \ne k(*)} b_m = c_m$ hence
$b_{k(*)},c_{k(*)} \notin A$.  For some $t < n(R)$ we have $\{b_m:m < \ell g
(\bar b),m \ne k(*)\}$ is disjoint to $A^{i(*)}_t$.  We can find a function
$\sigma$ from $\{0,\dotsc,\ell g(\bar b)-1\}$ to $\{0,\dotsc,n(R)-1\}$ such
that $\sigma(k) = t \equiv k = k(*),b_m \in A^{i(*)}_\ell \Rightarrow 
\sigma(m) = \ell$ and $\ell_1 \ne \ell_2 \and b_{\ell_1} \notin A \and 
b_{\ell_2} \notin A \Rightarrow \sigma(\ell_1) \ne \sigma(\ell_2)$.  
For $s \in \{1,2\}$ and
$r < n(R)$ let $B^s_r = A^{i(*)}_r \cup \{b_\ell:\sigma(\ell)=t\}$ if
$r \ne k(*)$ and $B^1_r = A^{i(*)}_r \cup \{b_{k(*)}\},B^2_r = A^{i(*)}_r
\cup \{c_{k(*)}\}$ if $r = k(*)$.  For each $s \in \{1,2\}$ we ask, 
choosing $\langle A^{i(*)+1}_r:r < n(R) \rangle$ as $\langle B^s_r:r < n(R) 
\rangle$ which of the demands hold.  Now $B^s_r$ extends $A^{i(*)}_r$ 
(so clause (v) holds), is a subset of ${\Cal U}$ with
$\le |A^{i(*)}_r| +1 \le i(*)+1$ element (by the choice of $\sigma$), so 
clauses (i) + (iii) holds and $r_1 \ne r_2 \Rightarrow B^s_{r_1} \cap
B^s_{r_2} = \emptyset$ (again look at the choice of $\sigma$) so clause (ii)
holds.  So necessarily clause (iv) fails.  For $r < n(R)$ let $E_r$ be the
following equivalence relation on $A^{i(*)}_r:a' E_r a''$ iff $a',a'' \in
A^{i(*)}_r$ and tp$_{\text{bs}}(a',\dbcu_{m \ne r} A^{i(*)}_m,R) =
\text{ tp}_{\text{bs}}(a'',\dbcu_{m \ne r} A^{i(*)}_m,R)$.  For $s \in
\{1,2\},k < n(R)$ let $E^s_r$ be the following equivalence relation on
$B^s_r:a' E_r a''$ iff $a',a'' \in B^s_r$ and tp$_{\text{bs}}(a',
\dbcu_{m \ne r} B^s_m,R) = \text{ tp}_{\text{bs}}(a'',\dbcu_{m \ne r}
B^s_m,R)$.
\sn
Now by the definition of $E_r$ clearly

$$
|\{\text{tp}_{\text{bs}}(a,\dbcu_{m \ne r} A^{i(*)}_m,R):a \in A^{i(*)}_r\}|
= |A^{i(*)}_r/E_r|
$$
\mn
hence as $\langle A^{i(*)}_r:r < n(R) \rangle$ satisfies $(i) - (iv)$ we
know that
\mr
\item "{$(*)_1$}"  $i(*) \le \dsize \sum_{r < n(R)} |A^{i(*)}_r/E_k|$.
\ermn
Also

$$
|\{\text{tp}_{\text{bs}}(a,\dbcu_{m \ne r} B^s_m,R):a \in B^s_r\}| = 
|B^s_r/E^s_r|
$$
\mn
hence as $\langle B^s_r:k < n(R) \rangle$ fail condition $(iv)$ (see
above) we have
\mr
\item "{$(*)_2$}"  $i(*)+1 > \dsize \sum_{r < n(R)} |B^s_r/E_r|$.
\ermn
Now for each $r < n(R)$, clearly $E^s_r \restriction A^{i(*)}_r$ is an
equivalence relation refining $E_r$, hence
\mr
\item "{$(*)_3$}"  $|A^{i(*)}_r/E_r| \le |A^{i(*)}_r/E^s_r| \le 
|B^s_r/E^s_r|$.
\ermn
The three together gives
\mr
\item "{$(*)_4$}"   $|A^{i(*)}_r/E_r| = |A^{i(*)}_r/E^s_r| = |B^s_r/E^s_r|$
\ermn
hence
\mr
\item "{$(*)_5$}"   $E_r = E^s_r \restriction A^{i(*)}_r$
\sn
\item "{$(*)_6$}"   if $d \in B^s_r \backslash A^{i(*)}_r$ then for some
$d' \in A^{i(*)}_r$ we have $d E^s_k d'$.
\ermn
Apply $(*)_6$ to $r=t$ choosing $d_s = b_{k(*)}$ if 
$s=1$ and choosing $d_s = c_{k(*)}$ if $s=2$ so
$d_s \in B^s_t \backslash A^{i(*)}_t$ hence there is $d'_s \in A^{i(*)}_t$ 
such that $d_s E^s_t d'_s$ so tp$_{\text{bs}}(d_s,\dbcu_{m \ne t} B^s_m,R) =
\text{tp}_{\text{bs}}(d'_s,\dbcu_{mr \ne t} B^s_m,R)$ 
hence tp$_{\text{bs}}(d_s,\dbcu_{m \ne t} A^{i(*)}_m,R) = 
\text{ tp}_{\text{bs}}(d'_s,\dbcu_{m \ne t} A^{i(*)}_m,R)$.  \nl
But tp$_{\text{bs}}(d_1,\dbcu_{m \ne t}
A^{i(*)}_m,R) = \text{tp}_{\text{bs}}(b_{k(*)},\dbcu_{m \ne t} A^{i(*)}_m,R) =
\text{ tp}_{\text{bs}}(c_{k(*)},\dbcu_{m \ne t} A^{i(*)}_m,R) =$ \nl
$\text{tp}_{\text{bs}}(d_2,\dbcu_{m \ne t} A^{i(*)}_t,R)$ (second equality
as $b_{k(*)} E_A c_{k(*)}$ by the choice of $\bar b,\bar c$).
\mn
So together with the previous sentences
$\text{tp}_{\text{bs}}(d'_1,\dbcu_{m \ne t} A^{i(*)}_m,R) =
\text{ tp}_{\text{bs}}(d'_2,\dbcu_{m \ne t} A^{i(*)}_m,R)$ that is
$d'_1 E_t d'_2$ (recall $d'_1,d'_2 \in A^{i(*)}_t$).  So by $(*)_5$ we have 
$d'_1 E^s_t d'_2$ for $s=1,2$.  Clearly $m < n(R) \and m \ne t \Rightarrow
B^1_m = B^2_m$ hence $\dbcu_{m \ne t} B^1_m = \dbcu_{m \ne t} B^2_m$
and let $E^*_t$ be the following equivalence relation on
${\Cal U}:a' E^*_t a''$
iff tp$_{\text{bs}}(a',\dbcu_{m \ne t} B^s_m,R) = \text{tp}_{\text{bs}}
(a'',\dbcu_{m \ne t} B^s_m,R)$.  Clearly $E^s_t = E^*_s \restriction B^s_t$,
hence $d_1 E^*_t d'_1$ (by the choice of $d'_1$), $d'_1 E^*_t d'_2$ (see the
previous sentences) and $d'_2 E^*_t d_2$ (by the choice of $d'_2$).  Together 
as $d_1 = b_{k(*)},d_2 = c_{k(*)}$ we have $b_{k(*)} E^*_t c_{k(*)}$; 
but $\{c_\ell:\ell \ne k(*)\} =
\{b_\ell:\ell \ne k(*)\} \subseteq \dbcu_{m \ne t} B^s_m$ (by the choice of
$\sigma$) so $\text{tp}_{\text{bs}}(b_{k(*)},\{c_\ell:\ell \ne k(*)\},R) =
\text{ tp}_{\text{bs}}(c_{k(*)},\{c_\ell:\ell \ne k(*)\},R)$ 
a contradiction to the choice of $\bar b,\bar c$.

So $A = \dbcu_{r < n(R)} A^{i(*)}_r$ satisfies clause (B) of \scite{3.4}. 
Note that $E_A$ has $\le |A| + \lambda_1(R)$ equivalence classes by the
definition of $\lambda_1(R)$, so $A$ satisfies clause (C), and $\bar C$ 
satisfies clause (D) so is really as required.  \hfill$\square_{\scite{3.4}}$
\enddemo
\bigskip

\demo{\stag{3.5} Conclusion}  Letting $\lambda_\ell = \lambda_\ell(R)$ we
have $\exists_R$ is bi-interpretable with \nl
$\{Q^{\text{mon}}_{\lambda_0},
Q^{1-1}_{\lambda_1},\exists_{R_1},\exists_E\}$, where $|\text{Dom}(R_1)| 
\le n(R)^2 \lambda_1(R)$
and $E$ is an equivalence relation on ${\Cal U}$.  This is done uniformly 
(i.e. the formulas depend on $n(R)$ only).
\enddemo
\bigskip

\remark{Remark}  Note that $Q^{\text{mon}}_{\lambda_0}$ can be omitted being
swallowed by $\exists_E$.
\endremark
\bigskip

\demo{Proof}  We've shown $Q^{1-1}_{\lambda_1(R)} \le_{\text{int}} 
\exists_R$ (see \scite{3.3}).  Let $A,C$ be as in the lemma \scite{3.4},
choose $A_1$ such that $A^1 \cap A = \emptyset,|A^1| \le n(R)^2 \lambda_1(R)$
and $A \cup A^1$ includes $\ge \text{Min} \{n(R),|a/E_A|\}$ elements of each 
$E_A$ 
equivalence class $a/E_A$.  Lastly let $R_1 = R \restriction (A \cup A^1)$.
\nl
Now by the choice of $A$ and $\bar C$ clearly

$$
R(x_1,\dotsc,s_{n(R)}) \text{ iff}
$$

$$
\align
(\exists y_1) \cdots (\exists y_{n(R)})(\dsize \bigwedge_{1 \le i \le n(R)}
x_i E_{A,\bar C} y_i &\and \dsize \bigwedge_{i,j=1,\dotsc,n(R)} x_1 = x_j
\equiv y_i = y_j \\
  &\qquad \qquad \qquad \quad \and R_1(\bar y)).
\endalign
$$
\mn
So $\exists_R \le_{\text{int}} \{Q^{1-1}_{\lambda_1},\exists_{R_1},
\exists_{E_{A,\bar C}}\}$.  Now
$\exists_{R_1} \le_{\text{int}} \{\exists_R,Q^{\text{mon}}_{\lambda_1}\}$
by the definition of $R_1,\exists_{E_{A,\bar C}} \le_{\text{int}} \{\exists_R,
Q^{\text{mon}}_{\lambda_1}\}$ directly and $Q^{1-1}_{\lambda_1} \le
\exists_R$ by \scite{3.3} and $Q^{\text{mon}}_{\lambda_1} \le Q^{1-1}
_{\lambda_1}$.  So
$\{Q^{1-1}_{\lambda_1},Q^{\text{mon}}_{\lambda_0},\exists_{R_1},
\exists_{E_A}\} \le_{\text{int}} \{\exists_R,Q^{\text{mon}}{\lambda_0}\}$ and
we finish.  \hfill$\square_{\sciteu{3.5}}$\sciteuphantom{3.5}
\enddemo
\bigskip

\remark{\stag{3.6} Remark}  Note the $Q^{1-1}_{|\text{Dom}(R_1)|}$ is 
uniformly interpretable (for fixed $n(R))$ in $Q^{1-1}_{\lambda_1}$ including
the case $\lambda_1$ is finite, so \sciteu{3.5} holds for it too.
\endremark
\bigskip

\proclaim{\stag{3.7} Claim}  If $|{\Cal U}| > \lambda_1 \ge \lambda^k,R$ 
a $k$-place relation on $A \subseteq {\Cal U}$ and $|A| \le \lambda$
(and ${\Cal U}$ finite) \ub{then} $Q_R \le_{\text{exp}} Q^{1-1}_{\lambda_1}$.
\endproclaim
\bigskip

\demo{Proof}  By \scite{1.12}.
\enddemo
\bigskip

\demo{\stag{3.8} Conclusion}  If $R$ is an 
$n(R)$-place relation on ${\Cal U}$ and
$\lambda_1(R)^{n(R)} \le |{\Cal U}|$, \ub{then} for some equivalence relation
$E$ we have

$$
\{Q^{\text{eq}}_E,Q^{1-1}_{\lambda_1(R)}\} \le_{\text{int}}
Q_R \le_{\text{exp}} \{Q^{\text{eq}}_E,Q^{1-1}_{\lambda_1(R)^{n(R)}}\}.
$$
\enddemo
\bigskip

\demo{Proof}  We have by \scite{3.3} that $Q^{1-1}_{\lambda_1(R)} \le
\exists_R$.  By \scite{3.7} for every binary relation $S$ on 
$\sqrt {\lambda_1}$ we have $\exists_S \le Q^{1-1}_{\lambda_1(R)}$.  
So every relation
on ${\Cal U}^{\frac{1}{2n(R)}}$ is interpreted in $\exists_R$.
\enddemo
\bigskip

\remark{\stag{3.9} Remark}  So up to 
expressability and up to a power by $n(R)$ (and
possibly increasing ${\Cal U}$), we have that
$\{Q^{\text{eq}}_E,Q^{1-1}_{\lambda_1(R)}\}$ exhaust all the
information on $Q_R$ (up to interpretability).
\endremark
\bn
We can get the parallel result for $Q_K$.
\definition{\stag{3.10} Definition}  $\lambda_1(K) = \{\lambda:
\text{ for every } R \in K \text{ we have } \lambda_1(R) < \lambda\}$.  Note
that the maximum is taken for each ${\Cal U}$ separately.
\enddefinition
\bigskip

\demo{\stag{3.12} Conclusion}  1) $Q^{1-1}_{< \lambda_1(K)} \le_{\text{int}}
\exists_K$. \nl
2) There are $K_1$ and $\bold E$, a family of equivalence relations (for
each ${\Cal U} \in {\frak U}$, closed under permutations of ${\Cal U}$) such
that:
\mr
\item "{$(a)$}"  $\exists_K \equiv_{\text{int}} \{\exists_{K_1},Q^{1-1}
_{< \lambda_1(K)},Q_{\bold E}\}$
\sn
\item "{$(b)$}"  for any $R \in K_1$ we have $|\text{Dom}(R)| < n(R)^2 \times
\mu(K)$ where $\mu(K) = \text{ Min}\{\mu:R \in K \Rightarrow |\text{Dom}(R)|
< \mu\}$ the minimum taken for each ${\Cal U} \in {\frak U}$ separately.
\endroster
\enddemo
\bigskip

\demo{Proof}  Straight by uniformity.
\enddemo
\newpage

\head {Assignments} \endhead  \resetall 
\bn
1) 2.10; 2 copies
\sn
2) 639a

\newpage
    
REFERENCES.  
\bibliographystyle{lit-plain}
\bibliography{lista,listb,listx,listf,liste}

\def\germ{\frak} \def\scr{\cal}
  \ifx\documentclass\undefinedcs\def\rm{\fam0\tenrm}\fi
  \def\defaultdefine#1#2{\expandafter\ifx\csname#1\endcsname\relax
  \expandafter\def\csname#1\endcsname{#2}\fi} \defaultdefine{Bbb}{\bf}
  \defaultdefine{frak}{\bf} \defaultdefine{mathbb}{\bf}
  \defaultdefine{beth}{BETH} \def\bbfI{{\Bbb I}} \def\mbox{\hbox}
  \def\text{\hbox} \def\om{\omega} \def\Cal#1{{\bf #1}} \def\pcf{pcf}
  \defaultdefine{cf}{cf} \defaultdefine{reals}{{\Bbb R}}
  \defaultdefine{real}{{\Bbb R}} \def\restriction{{|}} \def\club{CLUB}
  \def\w{\omega} \def\exist{\exists} \def\se{{\germ se}} \def\bb{{\bf b}}
  \def\equivalence{\equiv} \let\lt< \let\gt> \def\cite#1{[#1]}
  \def\implies{\Rightarrow}
\begin{thebibliography}{BlSh 156}
\makeatletter \renewcommand{\@biblabel}[1]{[#1]} \makeatother

\bibitem[Bl]{Bl}John~T. Baldwin.
\newblock {Definable second order quantifiers}.
\newblock In J.~Barwise and S.~Feferman, editors, {\em {Model Theoretic
  Logics}}, Perspectives in Mathematical Logic, chapter XII, pages 446--478.
  Springer-Verlag, New York Berlin Heidelberg Tokyo, 1985.

\bibitem[BlSh 156]{BlSh:156}John~T. Baldwin and Saharon Shelah.
\newblock {Second-order quantifiers and the complexity of theories}.
\newblock {\em {Notre Dame Journal of Formal Logic}}, {\bf 26}:229--303, 1985.
\newblock Proceedings of the 1980/1 Jerusalem Model Theory year.

\bibitem[Sh:e]{Sh:e}Saharon Shelah.
\newblock {\em {Non--structure theory}}, accepted.
\newblock {Oxford University Press}.

\bibitem[Sh 28]{Sh:28}Saharon Shelah.
\newblock {There are just four second-order quantifiers}.
\newblock {\em {Israel Journal of Mathematics}}, {\bf 15}:282--300, 1973.

\bibitem[Sh 171]{Sh:171}Saharon Shelah.
\newblock {Classifying of generalized quantifiers}.
\newblock In {\em {Around classification theory of models}}, volume 1182 of
  {\em {Lecture Notes in Mathematics}}, pages 1--46. {Springer, Berlin}, 1986.

\bibitem[Sh:F334]{Sh:F334}{Shelah, Saharon}.
\newblock {On quantification with a finite universe II}.

\end{thebibliography}

\enddocument

\bye